\documentstyle{amsppt}

\magnification=\magstep1
\hsize=14truecm
\vsize=22truecm

\define\Mr{M_n(\Bbb C)}
\define\Glin{Gl_n(\Bbb C)}
\define\ot{\otimes}
\define\ota{\otimes_\alpha}
\define\wten{\overset\vee\to\otimes}
\define\pten{\widehat\ot}
\define\iU{{U,V,\lambda}}
\define\FX{\Cal F(X)}
\define\FY{\Cal F(Y)}
\define\FZ{\Cal F(Z)}
\define\KX{\Cal K(X)}
\define\KY{\Cal K(Y)}
\define\KZ{\Cal K(Z)}
\define\A{$\Bbb A $}

\define\conv{\longrightarrow}
\define\subs{\subseteq}

\define\r{\rangle}
\redefine\span{\operatorname{span}}

\define\op{\text{op}}
\define\ginG{{g\in\Cal G}}
\define\bu{\bullet}
\define\dab{d^{\alpha ,\beta }}
\define\dia{d_{11}^\alpha}
\define\diab{d_{11}^{\alpha ,\beta }}
\define\da{d^\alpha}

\topmatter 
\title Amenability of Banach algebras\\ of\\ compact
operators \endtitle 
\author N\. Gr{\o}nb{\ae}k,\quad B\. E\.
Johnson,\quad and \quad G\. A\. Willis
\thanks GAW partially supported by SERC grant GR-F-74332.
\endthanks
\endauthor

\abstract
In this paper we study conditions on a Banach space $X$ that ensure
that the Banach algebra $\KX$ of compact operators is amenable.  We
give a symmetrized approximation property of $X$ which is proved to be
such a condition. This property is satisfied by a wide range of Banach
spaces including all the classical spaces. We then investigate which
constructions of new Banach spaces from old ones preserve the property
of carrying amenable algebras of compact operators. Roughly speaking,
dual spaces, predual spaces and certain tensor products do inherit
this property and direct sums do not. For direct sums this question is
closely related to factorization of linear operators. In the final
section we discuss some open questions, in particular, the converse
problem of what properties of $X$ are implied by the amenability of
$\KX$. 
\endabstract
\endtopmatter

\NoBlackBoxes

\heading
{0. Introduction }
\endheading
Amenability is a cohomological property of Banach algebras
which was introduced in [J]. The definition is given
below. It may be thought of as being, in some ways,
a weak finiteness condition. For example, amenability of
C*-algebras is equivalent to nuclearity, see \cite{Haa}. 
Also, a group algebra, $L^1(G),$ is amenable if and only if
the locally compact group, $G,$ is amenable, see
\cite{J}, and many theorems valid for finite or compact
groups have weaker generalizations to amenable groups but
to no larger class. This equivalence is the origin of the
term for Banach algebras. However, in some situations
amenability is not a finiteness condition. For example, a
uniform algebra is amenable if and only if it is
self-adjoint, see \cite{Sh}, and, for finite dimensional
Banach algebras, amenability is equivalent to
semisimplicity.  

The significance of amenability for some classes of Banach 
algebras suggests the question as to what it means for
other Banach algebras. In this paper we investigate the
amenability of the algebras of compact and of approximable
operators on the Banach space $X$. This was begun in [J],
where it is shown that ${\Cal K}(X)$ is amenable if $X$ is
$\ell_p,$ $1<p<\infty ,$ or $C[0,1].$ (${\Cal K}(X)$
denotes the algebra of compact operators on $X$ and ${\Cal
F}(X)$ the algebra of approximable operators.)  Relevant
properties of Banach spaces, such as the approximation
property, are now understood better than they were when
[J] was written and so we are able to make more progress.

We have not yet found such clear characterizations of
amenability for the algebras of approximable and compact 
operators as are known for classes of algebras mentioned
in the first paragraph. It does appear though that
amenability of ${\Cal F}(X)$ and ${\Cal K}(X)$ may be
equivalent to approximation properties for $X.$ One
immediate observation is that, since amenable Banach
algebras have bounded approximate identities, if the
algebra of compact operators on $X$ is amenable, then, by
\cite{D, Theorem 2.6}, $X$ has the bounded compact
approximation property and, if the algebra of approximable
operators is amenable, then $X$ has the bounded
approximation property. Moreover, results in \cite{G\&W}
and  \cite{Sa} show that, if ${\Cal K}(X)$ is
amenable, then $X^*$ has the bounded compact
approximation property and, if ${\Cal F}(X)$ is
amenable, then $X^*$ has the bounded approximation
property. It follows that, if ${\Cal F}(X)$ is amenable,
then ${\Cal K}(X) = {\Cal F}(X).$ 

Amenability of ${\Cal F}(X)$ is not equivalent to $X$ or $X^*$ having
the bounded approximation property however, as examples in the paper
show. Some sort of symmetry also seems to be required. In Section 3 we
formulate a symmetrized approximation property, called property (\A) ,
such that, if $X$ has property (\A) , then ${\Cal F}(X)$ is amenable.
This formulation is an abstract version of the argument used in [J].
We show that, if $X$ has a shrinking, subsymmetric basis, then it has
property (\A) and hence ${\Cal F}(X)$ is amenable. Many spaces which
do not have such a basis also have property (\A) .

The necessity of some sort of symmetry becomes apparent
when we consider the stability of the class of spaces $X$
such that ${\Cal F}(X)$ is amenable. Subject to some
restrictions, this class of spaces is closed under tensor
products and taking duals, as is shown in Sections 2 and 5.
However, it is not closed under direct sums or passing to
complemented subspaces, see Section 6. The results in
Sections 5 and 6 depend on some new stability properties
for amenable Banach algebras which we establish in those
sections.

Many questions remain to be answered before we understand fully the
connection, if any, between amenability of ${\Cal F}(X)$ or ${\Cal
K}(X)$ and approximation properties of $X.$ These questions are
discussed in the last section of the paper. We do not investigate
other homological properties of ${\Cal F}(X)$ and ${\Cal K}(X).$ One
other such property has been studied in \cite{Ly}.

We now give the definition of amenability for Banach algebras. It is
made in terms of Banach modules and derivations. Recall that, for a
Banach algebra ${\Cal A},$ a Banach space $X$ is a Banach ${\Cal
A}$-bimodule if $X$ is a ${\Cal A}$-bimodule and there is a constant
$K$ such that $||a.x||\leq K ||a||\,||x||$ and $||x.a||\leq K
||a||\,||x||$ for each $a$ in ${\Cal A}$ and $x$ in $X.$ If $X$ is a
Banach ${\Cal A}$-bimodule, then the dual space, $X^*,$ is a Banach
${\Cal A}$-bimodule with the actions defined by $\langle a.x^*,x \r =
\langle x^*,x.a \r$ and $\langle x^*.a,x\r = \langle x^*,a.x\r ,$ for
$a$ in ${\Cal A},$ $x$ in $X$ and $x^*$ in $X^*.$ A derivation into an
${\Cal A}$-bimodule $X$ is a linear map $D:{\Cal A}\to X$ such that
$D(ab) = a.D(b) + D(a).b,$ for all $a,b$ in ${\Cal A}.$ If $x$ belongs
to $X,$ then the map $a\mapsto a.x - x.a$ is a derivation into $X.$
Such derivations are called inner.
 
\proclaim{Definition 0.1} The Banach algebra ${\Cal A}$ is
amenable if, for every Banach ${\Cal A}$-bimodule $X,$ every
continuous derivation $D:{\Cal A}\to X^*$ is inner. 
\endproclaim
   
\noindent See \cite{J, Section 5}, or \cite{B\&D, Definition VI.2}.   

This definition will sometimes be used directly but we
will often use another characterization of amenability,
namely that ${\Cal A}$ is an amenable Banach algebra if
and only if ${\Cal A}\pten {\Cal A}$ has an approximate
diagonal. An approximate diagonal is a bounded net, $\{
d_\lambda \} _{\lambda \in \Lambda },$ in 
${\Cal A}\pten {\Cal A}$ such that 
$$\lim _{\lambda \to \infty }||a.d_\lambda - d_\lambda
.a|| = 0\ \ \ \  \hbox{ and } \ \ \ \lim_{\lambda \to
\infty } ||\pi (d_\lambda )a - a|| = 0, \ \ \ \  
(a\in {\Cal A}), $$
where $\pi $ denotes the product map 
${\Cal A}\pten {\Cal A}\to {\Cal A}$ and module
actions on ${\Cal A}\pten {\Cal A}$ are defined by
$a.(b\otimes c) = (ab)\otimes c$ and $(b\otimes c).a =
b\otimes (ca),$ for $a, b $ and $c$ in ${\Cal A}.$ If we define a
product on ${\Cal A}\pten {\Cal A}$ by $(a\ot b)(c\ot d)=ac\ot db,$
then the first of these conditons can also be stated as
$$
\lim_{\lambda \to \infty}\Vert (a\ot 1 - 1\ot a)d_\lambda\Vert =0\quad
(a\in {\Cal A}), 
$$ 
where $1$ is a formally adjoined unit. Approximate diagonals are
useful, for example, when we show that, if $X$ has property (\A), then
${\Cal F}(X)$ is amenable.

\heading{1. Notation}
\endheading

We begin by establishing notation. Throughout, $X$ and $Y$ will denote
(infinite dimensional) Banach spaces and $X^*$ the space of bounded
linear functionals on $X$ with its usual norm. Small letters $x$ etc\.
will denote elements in $X$, whereas $x^*$ etc\. will denote elements
in $X^*$.  We will consider the following classes of operators:

$$\align 
&F(X,Y)= \{\text{\rm finite rank operators  $X\to Y$} \}\\
&\Cal N (X,Y)= \{\text{\rm nuclear operators $X\to Y$}\}\\
&\Cal F(X,Y)= \text{\rm uniform closure of } F(X,Y) \\
&\phantom{\Cal F(X,Y)}= \{\text{\rm approximable operators $X\to Y$} \}\\
&\Cal K(X,Y)= \{\text{\rm compact operators $X\to Y$} \}\\
&\Cal I(X,Y)= \{\text{\rm integral operators $X\to Y$} \}\\
&\Cal B(X,Y)= \{\text{\rm bounded operators $X\to Y$} \}
\endalign
$$
We shall  write $F(X)$ instead of $F(X,X)$ etc.  

These are all two-sided operator ideals in $\Cal B(X,Y)$,and when
$X=Y$ they are, except $F(X)$, Banach algebras in their natural norms.
We refer the reader to any of \cite{D\&U}, \cite{Pie}, \cite{Pis} for
details.

Finite rank operators will, when convenient, be written as tensors,
that is, if $x^*_1,\dots ,x^*_n \in X^*$ and $y_1,\dots ,y_n \in Y$,
we shall denote the operator $x\to \sum x^*_i(x)y_i$ by $\sum y_i\ot x^*_i$.

If $S \in \Cal B(X,Y)$ we denote the adjoint map in $\Cal B(Y^*,X^*)$
by $S^a$, i.e. 
$$ 
\langle S(x),y^* \r\,=\,\langle x,S^a(y^*) \r \qquad (x \in X,
y^* \in Y^*) 
$$

\noindent
If $M\subseteq \Cal B(X,Y)$ we define $M^a \subseteq \Cal B(Y^*,X^*)$ by 
$$
M^a=\{T^a\vert T \in M \}
$$
This should not be confused with the notation for dual space. For
instance, if $X$ has Grothendieck's approximation property, then
$\Cal N(X)^*=\Cal B(X^*)$, whereas $\Cal N(X)^a$ is the set of
so-called $X$-nuclear operators on $X^*$.

We shall use the concepts left approximate identity, bounded left
approximate identity etc. in accordance with \cite{B\&D}.

\heading{2. Tensor products}
\endheading

It is of course important to be able to form new Banach 
spaces from old ones while preserving the property of
carrying amenable algebras of compact operators. The
first case to be considered is that of taking tensor
products because many important spaces can be viewed as
appropriate tensor products, for instance $L_p$-spaces
with values in a Banach space. We shall here investigate
whether amenability of $\KX$ and $\KY$ implies
amenability of $\KZ$, when $Z$ is the completion of $X\ot
Y$ in some crossnorm topology. An obvious approach to
this problem is to try to show that $\KX\ot \KY$ is a
dense subalgebra of $\KZ$ and then to deduce amenability
of $\KZ$ from that of $\KX\pten\KY$ by an appeal to
\cite{J, Corollary 5.5}. This program is considerably
easier to carry through if $X$ and $Y$ have the
approximation property. However, rather than making this
assumption, we prefer to work with approximable operators
instead of compact operators. The definition to follow
describes what is needed for above mentioned program to
work.

\proclaim{Definition 2.1} Let $X$ and $Y$ be Banach spaces and let
$\alpha$ be a crossnorm on $X\ot Y$. Denote the completion by $X\ota
Y$. We call $X\ota Y$ a {\it tight\/} tensor product of $X$ and $Y$,
if the following two conditions hold.
\roster
\item"(i)" There is $K>0$ so that for all $S\in\FX,\,T\in\FY$ the
operator on $X\ot Y$ given by 
$$
(S\ot T)x\ot y = Sx\ot Ty\quad(x\in X,\,y\in Y)
$$
has $\alpha$- operator norm not exceeding $K\Vert S\Vert\Vert T\Vert$.
\item"(ii)" $\operatorname{span} \{S\ot T \mid S\in \FX, T\in\FY\}$ is
dense in
$\Cal F(X\ota Y)$
\endroster
\endproclaim  

\noindent{\smc Remark:} The condition (i) is apparently weaker than
Grothendieck's $\ot$-norm condition \cite{Gr,Ch.1.3} in that it only
concerns finite rank operators on a tensor product into itself.

With this definition we have the obvious:

\proclaim{Theorem 2.2} Suppose that $\FX$ and $\FY$ are amenable and
that $X\ota Y$ is a tight tensor product. Then $\Cal F(X\ota Y)$ is
amenable. 
\endproclaim
\demo{Proof} \cite{J,Corollary 5.5}
\enddemo 

To apply this theorem we need to be able to recognize tight tensor
products. The following easy proposition is helpful. It shows that, as
usual when dealing with tensor products, it is important to be able to
identify $(X\ota Y)^*$. We shall view $(X\ota Y)^*$ as a subspace of
${\Cal{B}}(Y,X^*)$ (or equivalently of ${\Cal{B}}(X,Y^*)$). We give
${\Cal{B}}(Y,X^*)$ the canonical structure as a right Banach module
over $\FX$ and $\FY$, that is, the module actions are the restrictions
of the canonical actions of $\Cal B(X)$ and $\Cal B(Y)$.

 \proclaim{Proposition 2.3} Let $X$ and $Y$ be Banach spaces and let
$\alpha$ be a crossnorm on $X\ot Y$. Then $X\ota Y$ is a tight tensor
product if and only if the following two conditions hold:
\roster
 \item"(i)" $(X\ota Y)^* $ is a right Banach $\FX$- and $\FY$-
submodule of $\Cal B(Y,X^*)$ 
 \item"(ii)" $X^*\ot Y^*$ is norm dense in $(X\ota Y)^*$.
 \endroster
 \endproclaim
\demo{Proof} (i): Let $z=\sum x_i\ot y_i\in X\ot Y$, let $\Phi\in (X\ota
Y)^*$, and let $S\in \FX,\, T\in\FY$. Then 
 $$
 \align
\langle S\ot T\,z,\Phi\r&=\sum \langle Sx_i\ot Ty_i,\Phi\r\\
&=\sum\langle Sx_i ,\Phi Ty_i \r\\
&=\langle z,S^a\Phi T\r,
 \endalign
 $$
so that $S\ot
T$ is
$\alpha$-bounded
with $\Vert
S\ot T\Vert_\alpha\leq K\Vert S\Vert\Vert T\Vert $ if and only if the
Banach  
module properties hold with module constants $K_XK_Y\leq K$. 

(ii):  We shall use the identification $\FZ 
=Z\wten Z^*$.The canonical map $\FX\ot\FY\to{\Cal{F}}(X\ota Y)$
then becomes
 $$
(x\wten x^*)\ot(y\wten y^*)\to(x\ota y)\wten(x^*\ot y^*).
 $$

Using the injective property of $\wten$, it is now clear that the
image of $\FX\ot\FY$ is dense if and only if $X^*\ot Y^*$ is dense in
$(X\ota Y)^*$.
\enddemo

Recall that a crossnorm is called {\it reasonable} if the dual norm is
also a crossnorm. In this case tightness is particularly easy to
describe.

 \proclaim{Corollary 2.4} Suppose that $\alpha$ is a reasonable
crossnorm on $X\ot Y$ and that the module property 2.3.(i) holds.
Then $X\ota Y$ is tight if and only if 
 $$
(X\ota Y)^*=X^*\ot_{\alpha ^*} Y^*,
 $$
where $\alpha^*$ denotes the dual norm. 
 \endproclaim

With Proposition 2.3 at hand we can now give conditions for tightness
for some important tensor products. From 2.3.(ii) it is not surprising
that the Radon-Nikodym property enters the picture.

 \proclaim{Theorem 2.5} Let $X$ and $Y$ be Banach spaces, let $[0,1]$
be the unit interval, and let $(\Omega ,\Sigma ,\mu )$ be a
$\sigma$-finite measure space. Then  
\roster
 \item"(W)" The following are equivalent: 
  {\widestnumber\item{\quad(iii)}\roster
   \item"\quad(i)" $X\wten Y$ is tight for all $X$.
   \item"\quad(ii)" $C([0,1],Y)$ is a tight tensor product of $C[0,1]$ and
$Y$
   \item"\quad(iii)" $Y^*$ has RNP.\endroster} 
 \item"(P)" $X\pten Y$ is tight if and only if
${\Cal{F}}(Y,X^*)={\Cal{B}}(Y,X^*)$.
 \item"(M)" $L_p(\mu ,X),\,1\leq p<\infty$ is a tight tensor product of
$L_p(\mu )$ and $X$ if and only if $X^*$ has RNP with respect to $\mu$.
 \endroster
 \endproclaim
\demo{Proof} The identification of $(X\ota Y)^*$ with a subspace of
$\Cal B(Y,X^*)$ gives in the cases (W) and (P) $\Cal I(Y,X^*)$ and
$\Cal B(Y,X^*)$  respectively, so the module property 2.3.(i) is
obvious for these tensor products. Next, let
$S\in{\Cal{F}}(L_p(\mu ))\text{ and } T\in \FX$. From the
proof of Proposition 2.3.(i) it follows that it is enough
to show the submultiplicativity of the module norm for
$\Phi$ belonging to a norm determining subset of $L_p(\mu
,X)^*$. Let $\frac1{p}+\frac1{p'}=1$. Since $L_{p'}(\mu
,X^*)$ is isometrically embedded in $L_p(\mu ,X)^*$ and
since $L_p(\mu ,X)$ is isometrically embedded in
$L_{p'}(\mu ,X^*)^*$ it suffices look at $\Phi \in\Cal
B(L_p(\mu),X^*)$ coming from an element ${\bold{g}}\in
L_{p'}(\mu ,X^*)$.With the identifications being made we
have 
 $$
\Phi (f)=\int_\Omega f{\bold{g}}\,d\mu \quad(f\in L_p(\mu )).
 $$
 Then for $S\in{\Cal{B}}(L_p(\mu ))\text{ and }T\in\Cal
B(X)$
 $$
T^a\Phi S(f)=\int_\Omega S(f)T^a\circ {\bold{g}}\,d\mu .
 $$
 An appeal to the vector valued version of H\"olders
inequality, gives the desired norm inequality.

We now consider the statement 2.3.(ii) in our three cases. First we
look at (W). Since $(X\wten Y)^*={\Cal{I}}(Y,X^*)$ we are asking
whether the finite rank operators $X^*\ot Y^*$ are dense in $\Cal
I(Y,X^*)$ in the integral norm. The implication (i) $\Rightarrow$ (ii)
is obvious and (iii) $\Rightarrow$ (i) is valid because, if $Y^*$ has
RNP, then $\Cal I(Y,X^*)=\Cal N(Y,X^*)$ isometrically,
\cite{D\&U,Theorem VI.4.8, Corollary VIII.2.10}.  The implication (ii)
$\Rightarrow$ (iii) is true because, under the assumption (ii),
$F(C[0,1],Y^*)$ is dense in $\Cal I(C[0,1],Y^*)$. (We are here using
the symmetric r\^oles of $C[0,1]$ and $Y$.) By
\cite{D\&U,Theorem VI.3.12, Corollary VIII.2.10} every absolutely
summing operator $C[0,1]\to Y^*$ is nuclear since, by Lemma 2.8 below,
$\Cal I(C[0,1],Y^*) = \Cal N(C[0,1],Y^*)$, again using the symmetric
r\^oles of $C[0,1]$ and $Y$ and identifying $C[0,1]^*$ with $M[0,1]$,
the Banach space of Radon measures on the unit interval with the total
variation norm.. The RNP of $Y^*$ is now the content of \cite{D\&U,
Corollary VI.4.6}.

In the case (P) we just have to observe that $(X\pten
Y)^*={\Cal{B}}(Y,X^*)$ and $\operatorname{cl}(X^*\ot Y^*) =
{\Cal{F}}(Y,X^*)$.

Finally, as already noticed, $L_{p'}(\mu ,X^*)$ is isometrically
isomorphic to a subspace of $L_p(\mu ,X)^*$. As a consequence $L_p(\mu
,X)$ is tight if and only if $L_p(\mu ,X)^*=L_{p'}(\mu ,X^*)$. But
this is equivalent to $X^*$ having RNP with respect to $\mu$,
\cite{D\&U,Theorem IV.1.1}.
\enddemo   

>From a classical theorem by Pitt \cite{Pit} we get an
immediate consequence of Theorem 2.5.(P).

 \proclaim{Corollary 2.6} $\ell_p\pten\ell_q$ is tight if and only if
$\frac1p+\frac1q<1$.
 \endproclaim

\proclaim{Corollary 2.7} If $X^*$ has RNP and $\Cal F(X^*)$ is
amenable, then $\Cal F(\Cal F(X))$ is amenable.
\endproclaim
\demo{Proof} By Corollary 5.3 below, amenability of
$\Cal F(X^*)$ forces amenability of $\Cal F(X)$. The
identification $\Cal F(X) = X\wten X^*$ shows that $\Cal
F(X)$ is a tight tensor product of $X$ and $X^*$.
\enddemo

We have not been able to find the technical observation needed above
in the literature.

 \proclaim{Lemma 2.8} Let $M(K)$ be the Banach space of Radon
measures on a compact Hausdorff space with the total variation norm
and let $\Phi :Y\to M(K)$ be a finite rank operator. Then the
integral and nuclear norms of $\Phi$ coincide.
 \endproclaim
\demo{Proof} Since $M(K)$ is a ${\Cal{L}}_{1,1+\varepsilon}$-space
(cf\. \cite{ L\&P, Definition 3.1 }) for all
$\varepsilon>0$, there is a finite dimensional subspace
with $\operatorname{rg}\Phi\subs V $ and a projection
$P:M(K)\to V\text{ with }\Vert P\Vert\leq 1+\varepsilon $.
 Since $V$ is finite dimensional we have
${\Cal{N}}(Y,V)={\Cal{I}}(Y,V)$ isometrically. If 
$\iota:V\to M(K)$ is the inclusion map we get
$$
 \align
\Vert\Phi \Vert_{\text{nucl}}&=\Vert\iota P\Phi \Vert_{\text{nucl}}\\
&\leq\Vert P\Phi \Vert_{\text{nucl}}\\
&=\Vert P\Phi \Vert_{\text{int}}\\
&\leq(1+\varepsilon)\Vert\Phi \Vert_{\text{int}},
 \endalign
$$
so that $\Vert\Phi \Vert_{\text{nucl}}\leq\Vert\Phi
\Vert_{\text{int}}$. The reverse inequality is always true. 
\enddemo

\heading{3. Diagonals for $M_n(\Bbb C)$}
\endheading

For many Banach spaces $X$, in particular  the classical spaces, it is
possible to prove that $\KX$ is amenable as a consequence of a
uniform local structure of $X$, that is, as a consequence of a
property of finite dimensional subspaces. Before we set the scenario
in which this approach will work we shall take a closer look at
finite dimensional spaces. It is well known and easy to prove that
$M_n(\Bbb C)$ is amenable. In this section we shall view this in
terms of faithful irreducible representations of finite groups.
However, rather than speaking about faithful representations we shall
consider finite subgroups of $Gl_n(\Bbb C)$. Likewise, we shall
express irreducibility as a property of the embedding of the group
into $M_n(\Bbb C)$.

\proclaim{Lemma 3.1} Let $\goth D: \Cal G \to Gl_n(\Bbb C)$ be an
$n$-dimensional representation of a group $\Cal G$. Then $\goth D$ is
irreducible if and only if $\operatorname{span} \goth D(\Cal G) = M_n(\Bbb
C)$. 
\endproclaim

\demo{Proof} We extend the representation to the group algebra $\Bbb
C \Cal G$. Then the lemma is an easy consequence of Jacobson's density
theorem, \cite{B\&D, Theorem 24.10}.
\enddemo

Henceforth we shall deal with finite subgroups of $\Glin$ spanning
the whole of $\Mr$. These we shall call {\it irreducible $(n\times
n)$-matrix groups}. The connection of such with amenability of $\Mr$ is
described in the following proposition. The symbols $e_{ij}$ denote as
usual the matrix units.
 
 \proclaim{Proposition 3.2} Let $\Cal G$ be a finite irreducible
$(n\times n)$-matrix group. Then $\frac1{\vert\Cal G\vert} \sum_{g\in
\Cal G} g\otimes g^{-1}$ is equal to the canonical diagonal
$d_0=\frac1n\sum_{i,j=1}^n e_{ij}\otimes e_{ji}$ for $\Mr$. The
canonical diagonal $d_0$ is the only element of $\Mr\ot\Mr$ which is
simultaneously a diagonal for $\Mr$ and for the opposite algebra
$\Mr^{op}.$
 \endproclaim

\demo{Proof} Let $d=\frac1{\vert\Cal G\vert} \sum_{g\in
\Cal G} g\otimes g^{-1}$. That d is  a diagonal for
$\Mr$  means 

$$
\sum_\ginG ag\otimes g^{-1} = \sum_\ginG g\otimes g^{-1}a\quad(a\in
\Mr)\tag 3.1
$$
and
$$
\pi(d)=\operatorname{I}.\tag3.2
$$
Likewise, $d$ being a diagonal for $\Mr^{op}$ means

$$
\sum_\ginG ga\otimes g^{-1} = \sum_\ginG g\otimes
ag^{-1}\quad(a\in\Mr)\tag 3.3
$$
and 
$$
\pi_{op}(d)=\operatorname{I},\tag3.4
$$
where $\pi_{op}$ is the opposite multiplication $\pi_{op}(a\ot b)=ba.$

Since $\operatorname{span} {\Cal{G}}=\Mr$ it is enough to consider $a \in
{\Cal{G}}$ and then exploit linearity. We prove \thetag{3.3}:

$$
 \align
\sum_\ginG ga\otimes g^{-1} &=\sum_\ginG ga\otimes a(ga)^{-1}\\
&=\sum_{u\in {\Cal{G}}a} u\otimes au^{-1}
 \endalign
$$

Since ${\Cal{G}}a = {\Cal{G}}$, \thetag{3.3} follows. The identity
\thetag{3.1} is proved similarly, and \thetag{3.2} and \thetag{3.4} are
obvious. 

Simple computations with matrix units show that $d_0$ satisfies all of
\thetag{3.1},\dots,\thetag{3.4}. Now let $d=\sum_i
a_i\ot b_i$ be any element satisfying \thetag{3.2} and
\thetag{3.3} and write $d_0=\sum_j a'_j\ot b'_j$. Then
$$\align
d=\sum_i a_i\ot b_i&=\sum_{i,j} a_i\ot b_j'a_j' b_i\\
&=\sum_{i,j} a_ia_j'\ot b_j' b_i\\
&=\sum_{i,j} a_ib_ia_j'\ot b_j'\\
&=\sum_j a_j'\ot b_j'=d_0,
\endalign
$$
finishing the proof. Note that \thetag{3.1} and \thetag{3.4} follows
automatically from  \thetag{3.2} and \thetag{3.3}, since  $d_0$
satisfies \thetag{3.1} and \thetag{3.4}.

\enddemo 

(The use of the average $\frac1{\vert\Cal G\vert} \sum_{g\in \Cal G}
g\otimes g^{-1}$ probably dates back to the early days of
representation theory. It is a refinement of this which gives the
equivalence of amenability of group algebras and the existence of
invariant means, \cite{J, Theorem 2.5})

\subheading{Example 3.3} We shall several times in the sequel use
irreducible matrix groups of the following kind. Let $\Cal H$ be a
group of $(n\times n)$ permutation matrices corresponding to a
transitive subgroup of the symmetric group $S_n$. Then 
$$
\Cal G =\{ D(t){\pmb \sigma}\mid t\in\{\pm 1\}^n, {\pmb \sigma}\in\Cal
H \} 
$$
is an irreducible $(n\times n)$-matrix group. If $\Cal H = S_n$, then
$\Cal G$ is called {\it the monomial group of degree $n$}.

\heading{4. Amenability as a consequence of an approximation property}
\endheading

In this section we shall develop a method to lift
uniformly the diagonals of a matrix algebra to form an
approximate diagonal for $\FX$. The idea is illustrated
by the example $X=L_p(\mu )$. Locally $L_p(\mu )$ looks
like $\ell_p^n$ so we have \lq local\rq\  diagonals.
Furthermore, these diagonals are uniformly bounded (by
1). Using a direct limit argument we can form an
approximate diagonal for all of $\Cal F(L_p(\mu))$. 

This approach will work for all the classical spaces. The  definition
below is customised to make it work in a rather general situation. To
formulate it let us first look at a finite biorthogonal system
$\{(x_i,x_j^*)\, \vert\, x_i \in X;\, x_j^* \in X^*;\, i,j=1,\dots ,n \}$.
Using this system we may define a map $E:\Mr\to\FX$ by
 $$
E((a_{ij}))=\sum_{i,j} a_{ij}\, x_i\otimes x_j^* .
 $$
By biorthogonality, $E$ is an algebra homomorphism.  

 \proclaim{Definition 4.1} Let $X$ be a Banach space. We say that $X$
has property (\A) if there is a net of finite biorthogonal systems 
$$
\{(x_{i,\lambda},x_{j,\lambda}^*)\,\vert\, x_{i,\lambda} \in X;\,
x_{j,\lambda}^* \in X^*;\, i,j=1,\dots ,n_\lambda
\}\quad(\lambda\in\Lambda)
$$
and corresponding maps
 $$
E_\lambda :M_{n_\lambda }(\Bbb
C)\to\FX\quad(\lambda\in\Lambda)
 $$
such that with $P_\lambda = E_\lambda (I_{n_\lambda })$ the following
hold

\roster
 \item"\A (i)" $P_\lambda \conv 1_X$ strongly
 \item"\A (ii)" $P_\lambda^a \conv 1_{X^*}$ strongly
  \item"\A (iii)" For each $\lambda$ there is an irreducible
$(n_\lambda \times n_\lambda )$-matrix group
$\Cal G _ \lambda$ such that
$$
\sup \{\Vert E_\lambda(g)\Vert_{\op}\, \vert\, g \in \Cal G_\lambda,
\lambda
\in \Lambda\} < \infty. 
$$
 \endroster
 \endproclaim

We now show how to lift the diagonals of the matrix algebras to $\FX$. 

  \proclaim{Theorem 4.2} Suppose $X$ has property (\A). Then $\Cal
F(X)$ is amenable. 
 \endproclaim

\demo{Proof}
With notation as in the description of property (\A),
define the net $(d_\lambda )_{\lambda \in
\Lambda }$ in ${\Cal{F}}(X)\pten\Cal F(X)$ by
 $$
d_\lambda =\frac 1{|\Cal G_\lambda |}\sum_{g\in \Cal
G_\lambda }E_\lambda (g)\pten E_\lambda (g^{-1})\quad
(\lambda \in \Lambda ).
 $$
By assumption this is a bounded net. Observing that 
$\pi (d_\lambda )=P_\lambda$, we conclude by \A (i) that
$(\pi (d_\lambda ))_{\lambda \in \Lambda }$ is a bounded
left approximate identity for  $\FX$. 

Let $F\in \FX$. Then 
$$
\align
F.d_\lambda - d_\lambda .F&= (F-P_\lambda
FP_\lambda).d_\lambda - d_\lambda.(F-P_\lambda
FP_\lambda)\\ &\phantom{= }+P_\lambda
FP_\lambda.d_\lambda -d_\lambda.P_\lambda FP_\lambda\\
&=(F-P_\lambda FP_\lambda).d_\lambda -
d_\lambda.(F-P_\lambda FP_\lambda), \endalign
$$
since $\frac 1{|\Cal G_\lambda |}\sum_{g\in \Cal
G_\lambda }g\ot g^{-1}$ is a diagonal for
$M_{n_\lambda}(\Bbb C)$. By \A(ii) $(P_\lambda)$ is a
bounded right approximate identity for $\FX$, so that
$F.d_\lambda - d_\lambda .F \conv 0$

\enddemo

\noindent{\smc Remark 4.2}.a. The condition of biorthogonality in
property (\A) is stronger than necessary. The following asymptotic
trace condition suffices to establish amenability:
 $$
\frac1{n_\lambda } \sum_i^{n_\lambda }\langle x_{i,\lambda },x_{i,\lambda
}^*\r \conv 1 
$$
along $\Lambda$. With this condition replacing biorthogonality all
statements in this section about property (\A) remain valid. We have
made no use of this greater generality and so do not give the details
here. However, if $X$ has property (\A), then
\A(i) implies that $X$ is a $\pi$-space and so probably
there are spaces which satisfy the weaker condition but
not property (\A). (Note that apparently there are
spaces with the bounded approximation property which are
not $\pi$-spaces, see the introduction to [C\&K].)

\smallskip
\noindent{\smc Remark 4.2}.b Conditions \A(i) and \A(ii) together
imply that $X$ is what might be called a ``shrinking
$\pi_\lambda$-space", compare with the discussion in [G\&W]. Thus, if
$X$ has a basis and $P_n$ is the projection onto the span of the first
$n$ basis elements, then $(P_n)$ satisfying \A(i) and
\A(ii) implies that the basis is a shrinking basis. If, furthermore,
$(P_n)$
satisfies \A(iii) with the monomial group of degree $n$, then the basis
is a symmetric basis, see \cite{L\&T, Ch. 3a}. In this
case $X$ will have property (\A), see also Theorem 4.5
below.

\smallskip
Property (\A) is preserved for some natural  Banach spaces formed
from the original space, as set forth in the next two theorems. This
will enable us to establish amenability of $\FX$ for a
large class of Banach spaces, including all the classical
spaces.

 \proclaim{Theorem 4.3}
Let $X$ be a Banach space. If $X^*$ has property (\A), then $X$ has
property (\A).
\endproclaim
\demo{Proof}
Let
 $$
   \{(x_{i,\lambda }^*,x_{i,\lambda}^{**})\}\quad(\lambda\in\Lambda)  
 $$
be a net of biorthogonal systems satisfying the conditions of (\A) with
respect to $X^*$. Let $\Cal U\text{ and }\Cal V$ be the sets of all finite
dimensional subspaces of $X$ and $X^*$ respectively, and let $U\in\Cal
U$ and $V\in\Cal V$ be given. By means of the principle of local
reflexivity \cite{L\&T}, choose a linear map
 $$
S_\iU :\text{span}\,(\{x_{i,\lambda}^{**}\mid
i=1,\dots,n_\lambda \}\cup U)\longmapsto X
 $$
such that
\roster
 \item"(1)" $\Vert S_\iU\Vert\le 2$.
 \item"(2)" $S_\iU\vert_U=1_U$.
 \item"(3)"
$\langle S_\iU\big(x_{i,\lambda}^{**}\big), x^*\r=\langle x^*,
x_{i,\lambda}^{**}\r$ for all $x^*\in
\operatorname{span}(\{x_{i,\lambda}^*\}\cup V).$
\endroster       

\smallskip
We order $\Cal U$ and $\Cal V$ by containment and $\Cal U\times\Cal
V\times\Lambda$ by the product order. By construction 
$$
\{(S_\iU x^{**}_{i,\lambda },x^*_{j,\lambda
})\mid i,j=1,\dots ,n_\lambda \}\quad ((\iU)\in\Cal U\times\Cal V\times
\Lambda )
$$
is a net of finite biorthogonal systems. We denote the
corresponding lifts of matrix algebras by $E_\iU$ and the
corresponding projections by $P_\iU$. Then
 $$
P_\iU=S_\iU P_\lambda^a\iota_X\,,
 $$
where $\iota_X$ is the canonical inclusion of $X$ into $X^{**}$ and
$P_\lambda $'s are the property (\A) projections for $X^*$ . 
Clearly $\{P_\iU\}$ is a bounded set and for $x\in U$ 
 $$
\aligned
 \Vert P_\iU x-x\Vert &= \Vert
S_\iU\big(P_\lambda^ax)-x\Vert\\
&= \Vert S_\iU\big(P_\lambda^a x-x\big)\Vert\\
&\le  2\Vert P_\lambda^a x-x\Vert\,,
\endaligned
 $$
where the two last steps follow from (1) and (2) above. Hence \A(i) is
satisfied. Similarly for $x^*\in V$
 $$
\aligned
 P_\iU^a(x^*)&=\sum x_{i,\lambda}^*\otimes
S_\iU(x_{i,\lambda}^{**}\big)(x^*)\\ 
&=\sum \langle S_\iU\big(x_{i,\lambda}^{**}\big),
x^*\r x_{i,\lambda}^*\\
&=\sum \langle
x^*,x_{i,\lambda}^{**}\r x_{i,\lambda}^* = P_\lambda^a\big(x^*\big)\,,
\endaligned
 $$
using (3), so that \A(ii) is satisfied. The supremum in \A(iii) is
increased by at most a factor 2, using the same irreducible matrix
groups: $\Cal G_\iU = \Cal G_\lambda$.

We have thus found a net of finite biorthogonal systems which
satisfies the conditions needed for property (\A). 
  \enddemo
\smallskip
Property (\A) also behaves nicely with respect to tensor products:

 \proclaim{Theorem 4.4} Let $X$ and $Y$ be Banach spaces and let $Z$
be a tight tensor product of $X$ and $Y$. If $X$ and $Y$ have
property (\A), then so does $Z$.  
 \endproclaim
\demo{Proof} We write $Z=X\otimes_\alpha Y$. Let $({\Cal{O}}_\lambda
)_{\lambda \in\Lambda }\text{ and }({\Cal{R}}_\mu )_{\mu \in M }$ be
property (\A) nets of biorthogonal systems for $X$ and $Y$
respectively. We define the tensor product $({\Cal{O}}_\lambda
\ot{\Cal{R}}_\mu )_{(\lambda ,\mu )\in\Lambda \times M}$ to be the
product ordered net of biorthogonal systems for $X\ota Y$ given as 
 $$
{\Cal{O}}_\lambda \ot {\Cal{R}}_\mu = \{(x\ot y,x^*\ot y^*)\mid
(x,x^*)\in{\Cal{O}}_\lambda , \,(y,y^*)\in {\Cal{R}}_\mu \}.
 $$

Using the identification $\Mr\ot M_p(\Bbb C)=M_{np}(\Bbb C)$,
one checks easily that the property (\A) lifts
 $$
E_{(\lambda ,\mu)}:M_{n_\lambda n_\mu }(\Bbb C)\to \FZ
 $$
are nothing but $E_{(\lambda ,\mu)}=C_\lambda
\ot D_\mu $, where $C_\lambda\text{ and }D_\mu$ are the 
lifts belonging to $X$ and $Y$ respectively. Hence \A(i)
holds for $Z$ and, since $X^*\ot Y^*$ is dense in $Z^*,$
we also have \A(ii).  To obtain \A(iii) it suffices to
notice that, if ${\Cal{G}}\text{ and }{\Cal{H}}$ are
irreducible $(m\times m)\text{- and }(n\times n)$- matrix
groups, then ${\Cal{G}}\ot{\Cal{H}}$ is an irreducible
$(mn\times mn)$- matrix group.

\enddemo
\smallskip
We shall now give some concrete examples of spaces with property (\A).
The first is very much in the spirit of \cite{J, Proposition 6.1}.
Recall that a basis $(x_n)_{n\in\Bbb N}$ in a Banach space $X$ is
called {\it subsymmetric} if  $(x_n)_{n\in\Bbb N}$ is unconditional
and equivalent to the basis sequence $(x_{n_i})_{i\in\Bbb N}$ for
every increasing sequence $(n_i)_{i\in\Bbb N}$, see \cite {L\&T, Ch.
3.a} and \cite{Si, Ch. 21}. 

 \proclaim{Theorem 4.5} Suppose that $X$ has a subsymmetric and
shrinking basis. Then $X$ has property (\A).
 \endproclaim
\demo{Proof} Let $(x_n)_{n\in\Bbb N}$ be a subsymmetric and shrinking
basis and let $(x_n^*)_{n\in\Bbb N}$ be the associated sequence of
coordinate functionals. Then 
$$
\{(x_i,x^*_j)\mid i,j
=1,\dots,n\}\quad(n\in\Bbb N)
$$
is a sequence of finite biorthogonal systems satisfying 
the conditions of property (\A). \A(i) is immediate,
\A(ii) follows from the basis being shrinking. To prove
\A(iii) we shall use the following observations.

Since $(x_n)_{n\in\Bbb N}$ is unconditional, the family of operators
of the form  
$$
U\biggl(\sum_{n\in\Bbb N}^\infty a_n x_n\biggr)=\sum_{n\in \eta} s(n)
a_n x_n,\tag 4.1
$$
where $\eta\subs\Bbb N\text{ and }s\in\{\pm 1\}^{\Bbb N}$, is
uniformly bounded, say by $K>0$. The subsymmetry means that the
family of operators of the form
 $$
A_{(m_i)(n_i)}(x) = \sum_{i=1}^\infty x^*_{m_i}(x)x_{n_i}\tag 4.2$$
is uniformly bounded, say by $M>0$. Here $(m_i)_{i\in\Bbb N}\text{
and }(n_i)_{i\in\Bbb N}$ are two arbitrary increasing sequences of
integers.

Let ${\Cal{G}}_n$ be the subgroup of the monomial group of degree $n$
defined by the permutation matrix $\pmb\sigma$ corresponding to the cyclic
permutation $(12\cdots n)$, i\.e\. 
 $$
{\Cal{G}}_n = \{D(t){\pmb\sigma} ^k \mid t\in \{\pm 1\}^n,\,
k=0,\dots,n-1\}, 
 $$
By Lemma 3.1 $\Cal G_n$ is an irreducible $(n\times n)$-matrix group. 
We write elements in $X$ as sequences. Then for $g=D(t){\pmb\sigma
}^k$ we have: 
 $$
 \align
&E(g)(\xi _1,\xi _2,\dots,\xi _n,\dots) =\\
&(t(1)\xi _{n+1-k},\dots,t(k)\xi _n,t(k+1)\xi _1,\dots,t(n)\xi
_{n-k},0,\dots)
 \endalign
 $$
We see that $E(g)$ has the form $E(g)=A_1U_1+A_2U_2$
for appropriate choices of operators $U_i$ of type
\thetag{4.1} and $A_i$ of type \thetag{4.2}. Hence the
supremum in \A(iii) does not exceed $2KM$.  
\enddemo

 \proclaim{Corollary 4.6} Let $X$ be a reflexive Orlicz
sequence space or a reflexive Lorentz sequence space.
Then $X$ has property (\A) and so $\FX$ is amenable.
\endproclaim
\demo{Proof} See \cite{L\&T, Ch.3.a} for a discussion
showing that these spaces satisfy the hypotheses of
Theorem 4.5.  
\enddemo

We shall now give substance to the remark that the setup
of property (\A) is customized to deal with the classical
spaces.

 \proclaim{Theorem 4.7} Let $K$ be a compact Hausdorff space and let
$(\Omega ,\Sigma ,\mu )$ be a measure space. Then
$C(K)$ and $L_p(\mu )\ ,1\leq p\leq\infty,$ have property (\A).
\endproclaim

\demo{Proof}  Since $C(K)^*=L_1(\mu_K )$ for a suitable
measure space $(\Omega_K ,\Sigma_K ,\mu_K )$ and $L_\infty(\mu
)=C(K_\mu ) $ for a suitable compact space $K_\mu $, it follows from
Theorem 4.3 that it is enough to consider $L_p(\mu )$ for $1\leq
p<\infty$. We shall give the proof in detail in the case of a
probability space, cf. the remark below. Let $\Cal S$ be a finite
collection of disjoint measurable subsets of $\Omega$ whose union is
all of $\Omega$. As it is customary in integration theory we order
such dissections by ${\Cal{S}}_1\prec{\Cal{S}}_2$ if every set in
$\Cal S_1$ is a union of sets from $\Cal S_2$. We define the
biorthogonal systems in $L_p(\mu )\times L_{p'}(\mu )\
\frac1p+\frac1{p'}=1$ by
 $$
{\Cal{O}}_{\Cal{S}}=\{(\biggl(\frac1{\mu
(L)}\biggr)^\frac1p\chi_L,\biggl(\frac1{\mu
(M)}\biggr)^\frac1{p'}\chi_M)\mid  L,M \in {\Cal{S}}\},  
 $$
where $\chi_\bullet$'s denote indicator functions. It is now a routine
matter to verify property (\A). Let $P_{\Cal{S}}$ be the property (\A) 
projections. For an indicator function $\chi_M$ we have 
 $$
 \align
P_{\Cal{S}}(\chi _M)&=\chi _M\\
P_{\Cal{S}}^a(\chi _M)&=\chi _M
 \endalign
 $$
whenever $\{M\}\prec{\Cal{S}}$, so \A(i) and \A(ii) are immediate.
To prove \A(iii), consider $\Cal S=\{M_1,\dots,M_n\}$ and define
${\Cal{G}}_{\Cal{S}}$ to be the monomial group of degree $n$. Let
$g=D(t){\pmb\sigma }\text{ where }t\in \{\pm 1\}^n\text{ and
}{\pmb\sigma }$ is a permutation matrix. First notice that 
 $$
\Vert\sum_{i=1}^n a_i\biggl(\frac1{\mu (M_i)}\biggr)^\frac1p\chi
_{M_i}\Vert^p =
\sum_{i=1}^n \vert a_i\vert^p. 
 $$
 
Using this we get for an arbitrary $f\in L_p(\mu )$
 $$
 \align
\Vert
E_{\Cal{S}}(g)f\Vert^p&=\Vert\sum_{i=1}^n\biggl(t(i)\biggl(\frac1{\mu(M_i)}
\biggr)^\frac1{p'}\int_{M_i}f\,d\mu
\biggr)\biggl(\frac1{\mu (M_{\sigma 
(i)})}\biggr)^\frac1p\chi _{M_{\sigma (i)}} \Vert^p\\ 
&=\sum_{i=1}^n
\biggl(\frac 1{\mu(M_i)}\biggr)^\frac p{p'}\biggl(\vert\int_{M_i}f\,d\mu
\vert\biggr)^p\\ 
&\leq\sum_{i=1}^n
\biggl(\frac1{\mu(M_i)}\biggr)^\frac p{p'}\mu(M_i)^{\frac
p{p'}}\int_{M_i}\vert f\vert^p\,d\mu
\quad\text{(H\"older Inequality)}\\ &=\sum_{i=1}^n
\int_{M_i}\vert f\vert^p\,d\mu \\ &\leq \Vert f\Vert^p.
 \endalign
 $$
\enddemo\noindent 

{\smc Remark 4.7.}a. A proof of the general case can be given along
the same lines but with added minor technicalities. Alternatively, we
may reduce it to the special case.  We are interested only in
finite-dimensional subspaces. Functions in  such a subspace are supported on a
$\sigma$-finite measure space. The corresponding complemented
subspaces of $L_p(\mu)$ have projection constants uniformly bounded by
$1$ and are isometrically isomorphic to $L_p$-spaces of probability
measures.
\smallskip

Combining this with Theorem 2.5 and Theorem 4.4 we get a
large collection of Banach spaces carrying amenable
algebras.

\proclaim{Corollary 4.8} Suppose that $X$ has property
(\A). If $X^*$ has RNP, then $C(K,X)$ has property (\A).
If $X^*$ has RNP with respect to $\mu $, then $L_p(\mu
,X)$ has property (\A) for $1\leq p\leq\infty$.
\endproclaim 
\smallskip\noindent

\proclaim{Corollary 4.9} $\Cal F(\ell_p\pten\ell_q)$ is
amenable if and only if $\frac1p+\frac1q<1$.
\endproclaim

\demo{Proof} By Corollary 2.6 $\Cal F(\ell_p\pten\ell_q)$
is amenable for $\frac1p+\frac1q<1$. In \cite{A\&F} it is shown that,
if $r\leq s$, then $\Cal B(\ell_r,\ell_s) = (\ell _r\pten \ell
_{s'})^*$ contains a complemented copy of $\Cal B(\ell_2)$ and thus
fails the approximation property, \cite{Sz}. Hence, when
$\frac1p+\frac1q\geq 1$, then $\Cal F(\ell_p\pten\ell_q)$ does not
have a bounded right approximate identity and is consequently not
amenable.
\enddemo

Probably it is too much to hope that amenability of $\Cal F(X)$ is
equivalent to $X$ having property (\A). Since the approximate diagonal
stemming from property (\A) is obtained by means of lifts of the
canonical diagonals of matrix algebras, it will have the approximate
versions of the extra properties \thetag{3.3} and \thetag{3.4}. It seems
unlikely that such approximate diagonals should always exist, once
amenability of $\FX$ is established. In Section 6 we will see examples
of spaces for which $\Cal F(X)$ is amenable but for which we do not
know whether $X$ has property (\A) or even the weaker property
mentioned in Remark 4.2.a.

\heading{5. Dual spaces }
\endheading

We have seen that property (\A) \ passes from a dual Banach space to its
predual. The following stability property for amenability implies an
extension of this fact, namely, that if the algebra of approximable
operators on a dual space is amenable, then the algebra of
approximable operators on any predual of the space is amenable.  It
also implies a similar, but weaker, result for the algebra of compact
operators.

\proclaim{Theorem 5.1} Let ${\Cal A}$ be an amenable
Banach algebra and ${\Cal I}$ be a closed, left ideal in
${\Cal A}$ which has a bounded two-sided approximate
identity. Then ${\Cal I}$ is amenable. 
\endproclaim

\demo{Proof}It is convenient to define a new product on $\Cal
A\pten\Cal A$ by $(a\ot b)\bu(c\ot d)=ac\ot bd$, that is, $\bu$ is the
usual product on $\Cal A\pten \Cal A$. Let $(d_\alpha )_{\alpha \in A
}$ be an approximate diagonal for ${\Cal A}$, let $(e_{\beta })_{\beta
\in B }$ be a bounded two-sided approximate identity for ${\Cal I},$
and put
$$
p_{\alpha \beta \gamma } = d_\alpha\bu (e_\beta\ot
e_\gamma )\quad (\alpha \in A;\ \beta ,\gamma \in B). 
$$
Then, since ${\Cal I}$ is a left ideal and 
$(d_\alpha )_{\alpha \in A }$ and  
$(e_{\beta })_{\beta \in B }$ are bounded nets, 
$p_{\alpha \beta \gamma }$ belongs to a bounded
subset of ${\Cal I} \pten {\Cal I}.$ 
\par
For each $c$ in ${\Cal I}$ we have 
$$
\align
&\limsup_\gamma\Vert c.p_{\alpha \beta \gamma } - p_{\alpha \beta
\gamma }.c\Vert=\\
&\limsup_\gamma\Vert (c\ot 1)\bu  d_\alpha\bu
(e_\beta\ot e_\gamma) -  d_\alpha\bu (e_\beta\ot
e_\gamma )\bu (1\ot c)\Vert = \\
&\limsup_\gamma\Vert((c\ot 1)\bu d_\alpha-d_\alpha\bu
(1\ot c))\bu (e_\beta\ot e_\gamma )\Vert =\\
&\limsup_\gamma\Vert(c.d_\alpha-d_\alpha.c)\bu
(e_\beta\ot e_\gamma )\Vert, \endalign $$
using  $\lim_\gamma( e_\gamma c-ce_\gamma)= 0$

Since $(d_\alpha)_{\alpha\in A}$ is an approximate
diagonal and $(e_\beta)_{\beta\in B}$ is bounded we get from the
inequality 
$$
||(cd_\alpha-d_\alpha c)\bu e_\beta\ot e_\gamma||\leq
||(cd_\alpha-d_\alpha c)||\,|| e_\beta||\,|| e_\gamma||
$$
that
$$
\lim_\alpha\, \limsup_\beta\, \limsup_\gamma ||(cd_\alpha-d_\alpha c)\bu
e_\beta\ot e_\gamma||=0\,, 
$$
and so
$$
\lim_\alpha\, \limsup_\beta\, \limsup_\gamma
\Vert  c.p_{\alpha \beta \gamma } - p_{\alpha \beta \gamma }.c
\Vert = 0.
$$ 
\par
Furthermore,
$$
\align
\lim_\alpha \lim_\beta \lim_\gamma \pi 
(p_{\alpha \beta \gamma })c &=\lim_\alpha \lim_\beta \lim_\gamma \pi 
( d_\alpha\bu (e_\beta\ot e_\gamma))c\\
&=\lim_\alpha \lim_\beta \pi (d_\alpha\bu(e_\beta \ot
1))c\\ &=\lim_\alpha\pi(d_\alpha)c\\
&=c,
\endalign
$$
where the second and third equality follow from $(
e_\beta)_{\beta\in B}$ being a left approximate identity
for $\Cal I$ and $\Cal I$ being a left ideal, and the last
from $(d_\alpha)_{\alpha\in \Cal A}$ being an approximate
diagonal. It follows that we may choose a net from
$\{p_{\alpha \beta \gamma }\mid \alpha \in A;\ \beta
,\gamma \in B \} $ which is an approximate diagonal for 
${\Cal I}.$ Therefore ${\Cal I}$ is amenable. 
\enddemo

This theorem is an improvement on the last assertion in
Proposition 5.1 in [J]. It may also be shown, by a similar
argument but with $d_\alpha(e_\beta\ot e_\gamma)$ in place
of  $d_\alpha\bu (e_\beta\ot e_\gamma),$ that,
if ${\Cal A}$ is an amenable Banach algebra and ${\Cal I}$
is a two-sided ideal in ${\Cal A}$ with a bounded left
approximate identity, then ${\Cal I}$ is amenable.

\proclaim{Corollary 5.2} Let $X$ be a Banach space such
that ${\Cal K}(X^*)$ is amenable and ${\Cal K}(X)$ has a
bounded two-sided approximate identity. Then ${\Cal
K}(X)$ is amenable.
\endproclaim
\demo{Proof} ${\Cal K}(X)^a,$ which is anti-isomorphic to 
${\Cal K}(X),$ is a closed left ideal in  ${\Cal K}(X^*)$
and has a bounded two-sided approximate identity.
\enddemo

Example 4.3 in [G\&W] provides a Banach space, $X,$ such
that ${\Cal K}(X^*)$ has a bounded two-sided approximate 
identity but ${\Cal K}(X)$ does not. This example
suggests that the hypothesis that ${\Cal K}(X)$ has a
bounded two-sided approximate identity is necessary.
However, if $X$ has the approximation property it is not.

\proclaim{Corollary 5.3} Let $X$
be a Banach space such that ${\Cal F}(X^*)$ is amenable.
Then ${\Cal F}(X)$ is amenable.
\endproclaim
\demo{Proof} Since ${\Cal F}(X^*)$ has a bounded left
approximate identity, $X^*$ has the bounded approximation
property, by [D, Theorem 2.6]. Hence, by [G\&W, Theorem
3.3], ${\Cal F}(X)$  has a bounded two-sided approximate
identity.
\enddemo

It is an open question, which is discussed further in
Section 7, whether amenability of ${\Cal K}(X)$ implies
that $X$ has the approximation property.

The converse to Corollary 5.2 holds if ${\Cal K}(X^*)$
has a bounded two-sided approximate identity. This fact
will follow from another stability property of
amenability. 

\proclaim{Theorem 5.4} Let ${\Cal A}$ be a Banach
algebra  which has a bounded two-sided approximate
identity and let ${\Cal I}$ be a closed, left ideal in
${\Cal A}$ which is amenable and has a bounded left
approximate identity for ${\Cal A}.$ Then ${\Cal A}$ is
amenable. 
\endproclaim
\demo{Proof} 
By Proposition 1.8 in [J], it will suffice to check that
all derivations from ${\Cal A}$ into duals of essential 
${\Cal A}$-bimodules are inner. (An ${\Cal A}$-bimodule
$Y$ is {\it essential} if $Y = \span\{ a.y.b :
a,b\in {\Cal A}; y\in Y \} $, because, with the hypothesis of a
bounded approximate identity, this last space is closed.)
\par
Let $D:{\Cal A}\to Y^*$ be a derivation, where $Y$ is an 
essential ${\Cal A}$-bimodule. Since ${\Cal I}$ is
amenable, there is $y^*$ in $Y^*$ such that $Da = a.y^* -
y^*.a$ for every $a$ in ${\Cal I}.$  Then the map, $\delta
:{\Cal A}\to Y^*,$ defined by $\delta a =  a.y^* - y^*.a$
is an inner derivation from ${\Cal A}$ and so $D-\delta $
is a derivation from  ${\Cal A}$ whose restriction to
${\Cal I}$ is zero.  
\par
Now let $a$ and $b$ belong to ${\Cal A}$ and let
$(e_\lambda )_{\lambda \in \Lambda}$ be a bounded net in
${\Cal I}$ which is a left approximate identity for 
${\Cal A}.$ Then, since ${\Cal I}$ is a left ideal, 
$$
\align
0 &= \lim _{\lambda}\bigl(D-\delta
\bigr)(ae_\lambda ).b  \\
&= \lim _{\lambda}\bigl(D-\delta
\bigr)(a).e_\lambda b\\ 
&= \bigl(D-\delta \bigr)(a).b, 
\endalign
$$
where the two first identities are true because $ D-\delta$ is a
derivation which annihilates $ {\Cal I},$ and the third  because  
$(e_\lambda )_{\lambda \in \Lambda}$ is a left
approximate identity.

It follows that $\langle b.y , \bigl(D-\delta \bigr)(a) \r = 0$
for every $y$ in $Y$ and $a$ and $b$ in ${\Cal A}.$ Since
$Y$ is essential, $D = \delta $ and is thus inner. 
\enddemo

The next result may now be proved in a similar way to
Corollary 5.2.

\proclaim{Corollary 5.5} Let $X$ be a Banach space such
that ${\Cal K}(X)$ is amenable and ${\Cal K}(X^*)$ has a
bounded two-sided approximate identity. Then ${\Cal
K}(X^*)$ is amenable.
\endproclaim

The argument of Proposition 6.1 in [J] shows, without
change, that ${\Cal K}(c_0)$ is amenable. It follows from
this corollary and the fact that ${\Cal K}(\ell_1)$ has a
bounded two-sided approximate identity that 
${\Cal K}(\ell_1)$ is amenable. Proposition 6.1 in [J]
does not yield this fact about $\ell _1$ directly,
although we have shown it in Section 2 by
modifying the argument in [J] suitably.
\par \ 
\par
\noindent
{\bf Example 5.6. } The requirement in Corollary 5.5 that ${\Cal K}(X^*)$
have a bounded two-sided approximate identity is necessary.
Since $\ell _2$ is reflexive it has the RNP.
Hence, by Theorem 2.5, $\ell _2\wten \ell _2$ is a tight
tensor product and so, by Theorem 2.2, ${\Cal F}(\ell
_2\wten \ell _2)$ is amenable. Now $(\ell _2\wten \ell
_2)^*$ is isomorphic to $\ell _2\pten \ell _2$ and
$(\ell _2\wten \ell _2)^{**}$ to ${\Cal B}(\ell _2).$ 
Since  ${\Cal F}(\ell _2\wten \ell _2)$ has a bounded
two-sided approximate identity, ${\Cal F}(\ell _2\pten
\ell _2)$ has a bounded left approximate identity, see \cite{G\&W,
Theorem 3.3}.  However, ${\Cal B}(H)$ does not have the approximation
property (see \cite{Sz}) and so ${\Cal F}(\ell _2\pten \ell _2)$ does
not have a bounded right approximate identity. Therefore, ${\Cal
F}((\ell _2\wten \ell _2)^*)$ is not amenable.

\heading{ 6. Direct sums }
\endheading

In the following it is necessary to use the algebra of 
double multipliers on a Banach algebra ${\Cal A}.$ A
double multiplier on ${\Cal A}$ is a pair of
bounded operators, $(L,R),$ on ${\Cal A}$ which
commute and satisfy, for all $a$ and $b$ in ${\Cal A}:$
 $ L(ab)= L(a)b;$ $R(ab)=aR(b);$ and
$aL(b) = R(a)b.$ Denote the set of all double multipliers
on  ${\Cal A}$ by $M({\Cal A}).$ Then $M({\Cal A})$ is a
Banach space with the obvious norm and sum and becomes a
Banach algebra when equipped with the product
$(L_1,R_1)(L_2,R_2)=(L_1L_2,R_2R_1).$ If $T=(L,R)$ is a
double multiplier on ${\Cal A},$ then $L(a)$ will be
denoted by $Ta$ and  $R(a)$ by $aT.$ 
\par
Each element, $a,$ of ${\Cal A}$ determines a
double multiplier, $(L_a,R_a),$ where $L_a$ and $R_a$ are
respectively the operators on ${\Cal A}$ of left and
right multiplication by $a.$ Similarly, if ${\Cal A}$ is
embedded as an ideal in a Banach algebra ${\Cal B},$ then
each element of ${\Cal B}$ determines a double multiplier
on ${\Cal A}.$ Thus each operator on the Banach
space $X$ determines a double multiplier on ${\Cal
K}(X)$ and on ${\Cal F}(X).$ Note also that there is
always an identity, $I,$  in $M({\Cal A}).$ 

Now let $P_1$ be an idempotent in $M({\Cal A})$ and put
$P_2 =I-P_1$ and ${\Cal A}_{ij}=P_i{\Cal A}P_j, $  $i,j =
1,2.$ Next put 
${\Cal A}_{11}^\circ = \pi ({\Cal A}_{12}\pten {\Cal
A}_{21})$ and ${\Cal A}_{22}^\circ  = \pi ({\Cal
A}_{21}\pten  {\Cal A}_{12}),$ where $\pi $ denotes the
product in ${\Cal A}.$ Then ${\Cal A}_{ii}^\circ $ is
isomorphic, as a linear space, to the quotient of ${\Cal
A}_{ij}\pten  {\Cal A}_{ji},$ $j\ne i,$ by ker$(\pi )\cap 
({\Cal A}_{ij}\pten  {\Cal A}_{ji}).$ Let $\| .\| ^\circ $
denote the quotient norm on ${\Cal A}_{ii}^\circ .$ 

In this section we prove a couple of abstract results
about the stability of amenability when ${\Cal A}$ is
cut down to ${\Cal A}_{11}$ by an idempotent in 
$M({\Cal A})$ and then apply them to the
case where ${\Cal A}={\Cal K}(X)$ for some Banach space
$X$ and $P_1$ is determined by a projection on $X.$ We
will thus establish some stability properties of
amenability of ${\Cal K}(X)$ under direct sums of Banach
spaces.

\proclaim{ Proposition 6.1 } Let ${\Cal A}$ and 
${\Cal A}_{ij},\ i,j = 1,2,$ be as above. Then ${\Cal A}$
has a bounded two-sided approximate identity if and only
if ${\Cal A}_{11}$ and ${\Cal A}_{22}$ have bounded
two-sided approximate identities and ${\Cal A}_{ij}$ is
an essential left ${\Cal A}_{ii}$- and right ${\Cal
A}_{jj}$-module, $i,j = 1,2.$
\endproclaim

\demo{Proof}
Let $\{ e_\lambda \} _{\lambda \in \Lambda }$ be a bounded
net in ${\Cal A}.$ Then $\{ P_1e_\lambda P_1 + 
P_2e_\lambda P_2 \} _{\lambda \in \Lambda }$ is a
two-sided approximate identity if and only if $\{ 
P_ie_\lambda P_i \} _{\lambda \in \Lambda }$ is a
two-sided approximate identity in ${\Cal A}_{ii},$ a left
approximate identity for ${\Cal A}_{ij}$ and a right
approximate identity for ${\Cal A}_{ji},$  $i = 1,2;\
j\ne i.$
\enddemo 

The first of the abstract results is the following 

\proclaim{ Theorem 6.2} Let ${\Cal A}$ and  
${\Cal A}_{ij},\ i,j = 1,2,$ be as above and suppose that 
${\Cal A}$ has a bounded two-sided approximate identity and
that  ${\Cal A}_{22} = {\Cal A}_{22}^\circ .$ Then ${\Cal
A}$ is amenable if and only if ${\Cal A}_{11}$ is amenable.
\endproclaim

\demo{Proof}
  The inclusion map ${\Cal A}_{22}^\circ \to {\Cal A}_{22}$
is continuous and is also a surjection. Hence, by the
open mapping theorem, $\| .\| ^\circ $ is equivalent to
the given norm on ${\Cal A}_{22}.$ Furthermore, since
${\Cal A}$ has a bounded two-sided approximate identity,
Proposition 6.1 shows that ${\Cal A}_{22}$ also has a
bounded two-sided approximate identity. Therefore there is
a bounded net $\{ c^\beta \}_{\beta \in B }$ in ${\Cal
A}_{21}\pten {\Cal A}_{12}$ such that  $\{ \pi (c^\beta
)\} _{\beta \in B }$ is a bounded approximate identity for ${\Cal
A}_{22}.$ The elements of this net have the form  $c^\beta
= \sum _i r_i^\beta \ot s_i^\beta .$ 
\par
Now suppose that ${\Cal A}_{11}$ is amenable and let $\{
d_{11}^\alpha \} _{\alpha \in A }$ be an approximate
diagonal for ${\Cal A}_{11}.$  We will show that ${\Cal A}$
is amenable by showing that it has an approximate diagonal
consisting of elements of  the form  
$$
d^{\alpha ,\beta } = d_{11}^\alpha  + c_\beta d_{11}^\alpha.
$$
Here we have equipped $\Cal A\pten\Cal A$ with the product $(a\ot
b)(c\ot d)=ac\ot db$ as 
described in the introduction. Note first of all that the set of all
such elements is bounded because $\| d^{\alpha ,\beta }\|
\leq \| d_{11}^{\alpha }\| (1 + \| c^{\beta }\| ).$ In order to prove
that an approximate diagonal can be constructed, we shall use the
following
$$
\aligned
&\lim_\alpha[(\pi(c)\ot 1)(a_{21}\ot 1) - (1\ot a_{21})c]\dia =\\
&\lim_\alpha[(a_{12}\ot 1)c -(1\ot \pi(c))(1\ot a_{12})]\dia =0,
\endaligned\tag 6.1
$$
for each $c\in\Cal A_{21}\pten\Cal A_{12}$ and $a_{ij}\in\Cal A_{ij}$.
It is enough to prove \thetag{6.1}
for $c$ an elementary tensor $b_{21}\ot b_{12}$. Then the first
expression equals $(b_{21}\ot 1)(b_{12}a_{21}\ot 1 - 1\ot
b_{12}a_{21})\dia$, which tends to zero, because $(\dia)$ is an
approximate diagonal for $\Cal A_{11}$. The other limit is obtained
analogously. 

We will show that 
$$
\lim _\beta \lim _\alpha \pi (d^{\alpha ,\beta })a = a\quad (a\in
{\Cal A}), \tag 6.2
$$
and 
$$
\lim _\beta\, \limsup _\alpha \Vert(a\ot 1 - 1\ot a)d^{\alpha ,\beta
}\Vert  = 0\quad (a\in \Cal A). \tag 6.3 
$$
This will imply that an approximate diagonal can be constructed from the 
$d^{\alpha ,\beta }$'s.

First we prove \thetag{6.2}.  If $a$ is in $\Cal A,$
then $a=P_1a + P_2a$ and \thetag{6.2} follows because we
have 
$$ \lim_\alpha \pi(\dab)P_1a
=\lim_\alpha\pi(\dia)P_1a =P_1a, $$
since $\pi(c_\beta\dia)\in \Cal A_{21}\pi(\dia)\Cal A_{12}$ and
$\pi(\dia)$ is a bounded approximate identity for $\Cal A_{11}$.  
Likewise
$$
\align
\lim_\beta\lim_\alpha \pi(\dab)P_2a&=\lim_\beta\lim_\alpha
\pi(c_\beta\dia)P_2a\\
& = \lim_\beta \pi(c_\beta)P_2a = P_2a, 
\endalign
$$
again since $\pi(c_\beta\dia)\in \Cal A_{21}\pi(\dia)\Cal A_{12}$.

Now we prove \thetag{6.3}. Since
$a=a_{11}+a_{12}+a_{21}+a_{22},$ where $a_{ij}$ is in
${\Cal A}_{ij},$ we may treat these  terms separately.
We have 
  $$ \align
&(a_{11}\ot 1 - 1\ot a_{11})d^{\alpha,\beta}=(a_{11}\ot 1 - 1\ot
a_{11})d^\alpha_{11}\\ 
&(a_{12}\ot 1 - 1\ot a_{12})d^{\alpha,\beta}=(a_{12}\ot 1)c_\beta
d^\alpha_{11} - (1\ot a_{12})d^\alpha_{11}\\
&(a_{21}\ot 1 - 1\ot a_{21})d^{\alpha,\beta}=(a_{21}\ot 1)
d^\alpha_{11} - (1\ot a_{21})c_\beta d^\alpha_{11}\\
&(a_{22}\ot 1 - 1\ot a_{22})d^{\alpha,\beta}=(a_{22}\ot 1 - 1\ot
a_{22})c_\beta d^{\alpha,\beta}
\endalign
$$

Clearly the first term tends to 0 as
$\alpha\conv\infty$. The second term may be rewritten as
$$
\bigl((a_{12}\ot 1)c_\beta -(1\ot\pi(c_\beta))(1\ot
a_{12})\bigr)d^\alpha_{11} +1\ot(a_{12}\pi(c_\beta)-a_{12})d^\alpha_{11} 
$$
so that, using \thetag{6.1} and that $\bigl(\pi(c_\beta)\bigr)_{\beta
\in B}$ is a bounded right approximate identity for $\Cal A_{12}$, the
statement \thetag{6.3} is true in this case.

The third term may be rewritten as 
$$
\bigl((\pi(c_\beta)\ot 1)(a_{21}\ot 1) - (1\ot
a_{21})c_\beta\bigr)d^\alpha_{11} +
\bigl((a_{21}-\pi(c_\beta)a_{21})\ot 1\bigr)d^\alpha_{11}
$$
and treated analogously.

For the fourth term it is enough to look at elements of the form
$a_{22}=b_{21}b_{12}$ since by assumption these elements span a dense
subset of $\Cal A_{22}$ and we are working with bounded nets. We then get
$$
\align
&(a_{22}\ot 1 - 1\ot a_{22})d^{\alpha,\beta}=\\
&(b_{21}\ot 1)(b_{12}\ot 1 - 1\ot b_{12})d^{\alpha,\beta} + (1\ot
b_{12}) (b_{21}\ot 1 - 1\ot b_{21})d^{\alpha,\beta},
\endalign
$$
so that this case follows from the two previous cases.

To prove the converse, suppose now that ${\Cal A}$ is
amenable and let $\{ d^\alpha \} _{\alpha \in A}$ be an
approximate diagonal for ${\Cal A}$. Using the multiplier
multiplication $(P_i\ot P_j)(a\ot b)=P_ia\ot bP_j$ and $(a\ot
b)(P_i\ot P_j)= aP_i\ot P_jb$, we define
$$
\diab= (P_1\ot P_1)\da(P_1\ot P_1 +c_\beta)\quad (\alpha \in A\,,
\beta \in B).
$$

First note that for an elementary tensor we have
$$
\align
\lim_\beta\pi((a\ot b)c_\beta)&=\lim_\beta\pi((a\ot b)(P_2\ot P_2)
c_\beta)\\
&=\lim_\beta aP_2\pi(c_\beta)P_2a\\
&=aP_2b\\
&=\pi((a\ot b)( P_2\ot P_2))
\endalign
$$
so that 
$$
\align
\lim_\beta\pi(\diab)&= P_1\lim_\beta \pi(\da(P_1\ot P_1)
+c_\beta)P_1\\
&=P_1\pi(\da(P_1\ot P_1 + P_2\ot P_2))P_1\\
&=P_1\pi(\da)P_1,
\endalign
$$
which is a bounded left approximate identity for $\Cal A_{11}$,
directed over $\alpha\in A$.

For $a_{11}\in \Cal A_{11}$ we have
$$
\align
(a_{11}\ot 1 - 1\ot a_{11})\diab&=(a_{11}\ot 1 - 1\ot a_{11})(P_1\ot
P_1)\da(P_1\ot P_1 +c_\beta)\\
&=(P_1\ot P_1)(a_{11}\ot 1 - 1\ot a_{11})\da(P_1\ot P_1 +c_\beta),
\endalign
$$
which tends to $0$ as $\alpha\to\infty$, because $(\da)_{\alpha\in
A}$ is an approximate diagonal for $\Cal A$. This concludes the proof
of the theorem.

\enddemo

We now give some applications of this theorem in the case
when ${\Cal A} = {\Cal K}(X).$ 

\proclaim{Theorem 6.3} Let $X$ be a Banach space. Then 
${\Cal K}(X)$ is amenable if and only if ${\Cal
K}(X\oplus {\Bbb C})$  is amenable.
\endproclaim

\demo{Proof } Let ${\Cal A}= {\Cal K}(X\oplus {\Bbb C})$
and $P_1$ be the projection of $X\oplus {\Bbb C}$ onto
$X$ with kernel ${\Bbb C}.$ Then $P_2 = I-P_1$ is the
rank one projection onto ${\Bbb C}$ with kernel $X.$
Hence ${\Cal A}_{22} = P_2{\Cal K}(X\oplus {\Bbb C})P_2$
is the one-dimensional algebra spanned by $P_2.$  

It is easily seen that ${\Cal A}_{22}={\Cal A}_{22}^\circ
.$ Furthermore, since ${\Cal A}_{22}$ has an identity,
${\Cal A}$ has a bounded two-sided approximate identity
if either ${\Cal A}$ or ${\Cal A}_{11}$ is amenable.
Theorem 6.2 now applies.

\enddemo

Many of the classical Banach spaces are isomorphic to
their direct sum with the one-dimensional space and 
are also isomorphic to their hyperplanes. For some time it
was an unsolved problem, the so-called `hyperplane
problem', whether every Banach space has this property.
However, it is now known (\cite{G\&M}) that there is a Banach space
which is not isomorphic to any proper  subspace and so the
above theorem has some content. 
\par
An important class of Banach spaces is the class of
${\Cal L}_p$-spaces, where $1\leq p\leq \infty ,$
which were introduced in [L\&P]. The Banach space $X$ is
said to be an ${\Cal L}_{p,\lambda}$-space if there is a constant
$\lambda >0$ such that for every finite dimensional
subspace, $B,$ of $X$ there is a finite dimensional
subspace, $C,$ of $X$ such that $B\subseteq C$ and
$d(C,\ell _p^n)\leq \lambda ,$ where $n=\hbox{ dim }C.$
(If $Y$ and $Z$ are isomorphic Banach spaces,
then $d(Y,Z)$ is inf$(\| T\| ,\| T^{-1}\| ),$ where the
infimum is over all invertible operators, $T,$ from $Y$
onto $Z.$) Some examples of ${\Cal L}_p$-spaces  are
$\ell _p$ and $L_p(0,1).$ We have already seen in Theorem
4.7 that the algebras of compact operators on these 
examples are amenable.

Theorem III(c) in [L\&R] shows that ${\Cal L}_p$-spaces satisfy
stronger conditions than they are defined to have. Thus, if $X$ is an
${\Cal L}_p$-space, then there is a constant $\lambda '>0$ such that
for every finite dimensional subspace, $B,$ of $X$ there are a finite
dimensional subspace, $C,$ of $X$ and a projection, $P,$ of $X$ onto
$C$ such that: $B\subseteq C,$ $d(C,\ell _p^n)\leq \lambda ',$ where
$n=\hbox{ dim }C;$ and $\| P\|<\lambda '.$ It follows that every
${\Cal L}_p$-space has the approximation property, and so ${\Cal
K}(X)={\Cal F}(X)$ whenever $X$ is an ${\Cal L}_p$-space. It follows
also that, if $X$ and $Y$ are infinite dimensional ${\Cal
L}_p$-spaces, then every $T$ in ${\Cal F}(X)$ is a product $T=UV,$
where $U:X\to Y$ and $V:Y\to X$ are compact operators. Furthermore, if
$X$ is an ${\Cal L}_p$-space, then $X^*$ is an ${\Cal L}_q$-space,
where $q^{-1}+p^{-1}=1$ (\cite{L\&R, Theorem III(a)}).  Hence $X^*$ has
the bounded approximation property and so ${\Cal F}(X)$ has a bounded
two-sided approximate identity (\cite{G\&W, Theorem 3.3}). We are now
ready to prove 

\proclaim{Theorem 6.4} Let $1\leq p\leq\infty$ and let $X$ be an
${\Cal L}_p$-space. Then ${\Cal F}(X)$ is amenable. 
\endproclaim

\demo{ Proof } Let ${\Cal A}={\Cal F}(\ell _p\oplus
X)$ and let $P_1$ be the the idempotent in $M({\Cal A})$
determined by the projection onto $\ell _p$ with kernel
$X.$ Then ${\Cal A}$ has a bounded two-sided approximate
identity because ${\Cal A}_{11}={\Cal F}(\ell _p)$ and
${\Cal A}_{22}={\Cal F}(X)$ do. Also, since each compact
operator on $X$ factors through $\ell_p,$ ${\Cal
A}_{22}^\circ ={\Cal A}_{22}.$ Therefore, since ${\Cal
F}(\ell_p)$ is amenable, ${\Cal F}(\ell _p\oplus X)$ is
amenable by Theorem 6.2. 

That ${\Cal F}(X)$ is amenable now follows from another
application of Theorem 6.2 because ${\Cal F}(\ell _p\oplus
X)$ has a bounded two-sided approximate identity and
every compact operator on $\ell_p$ factors through $X.$ 
\enddemo

The finite rank projections on ${\Cal L}_p$-spaces which
were described above almost show that these spaces have
property (\A). The projections may be used to produce a
net of biorthogonal systems satisfying \A (i) and \A
(iii). However, it is not clear that the net will satisfy
\A (ii). If it could be shown that ${\Cal L}_p$-spaces in
fact have property (\A), then there would be a direct
proof of the amenability of ${\Cal F}(X)$ for these
spaces. It seems that indirect arguments are needed to
establish many of the properties of ${\Cal L}_p$-spaces,
see the remark after the statement of Theorem III in
[L\&R], and so it may be that they do not have property
(\A) .  Some specific examples for which this may be
tested are the spaces $\ell_2\oplus \ell_p.$ For
$1<p<\infty ,$ $\ell_2\oplus \ell_p$ is a ${\Cal
L}_p$-space, see [L\&P], example 8.2, but it is not clear
that it has property (\A).

Theorem 6.2 may be used to show that ${\Cal F}(X)$ is
amenable for some other spaces which may fail to have
property (\A). Let $\{ n_k\} _{k=1}^\infty $ be a sequence
of positive integers and choose $p$ and $q$ with $1\leq
p,q <\infty $ and  $p\ne q.$ Put $X= (\oplus _{k=1}^\infty
\ell _p^{n_k} )_{\ell _q}.$ Then $X$ has the bounded
approximation property and so ${\Cal K}(X) = {\Cal
F}(X).$ If $\{ n_k\} _{k=1}^\infty $ is bounded, then $X$
is isomorphic to $\ell _q$ and so we will suppose that 
$\{ n_k\} _{k=1}^\infty $ is not bounded. Clearly
$X$ is isomorphic to a complemented subspace of 
$(\oplus _{k=1}^\infty\ell _p )_{\ell _q}$ and so every
$T$ in ${\Cal F}(X)$ is a product $T=UV,$ where $U$ is in 
${\Cal F}(X,(\oplus _{k=1}^\infty\ell _p )_{\ell _q})$ and
$V$ in ${\Cal F}((\oplus _{k=1}^\infty\ell _p )_{\ell _q},
X).$ 

We also have that every $T$ in ${\Cal F}(
(\oplus _{k=1}^\infty\ell _p )_{\ell _q})$ is a product
$T=UV,$ where $U$ is in 
${\Cal F}((\oplus _{k=1}^\infty\ell _p )_{\ell _q}, X)$
and $V$ in 
${\Cal F}(X,(\oplus _{k=1}^\infty\ell _p )_{\ell _q}).$ To
see this, for each $r$ let $P_r$ be the natural rank $r^2$
projection of $(\oplus _{k=1}^\infty \ell _p)_{\ell _q}$
onto $(\oplus _{k=1}^r \ell _p^r)_{\ell _q}.$ Then $\{
P_r\} _{r=1}^\infty $ is a bounded left approximate
identity  for ${\Cal F}( (\oplus _{k=1}^\infty\ell _p
)_{\ell _q})$ and so $T= \sum _{r=1}^\infty P_r T_r,$
where $\sum _{r=1}^\infty ||T _r|| < \infty .$ Since 
 $\{ n_k\} _{k=1}^\infty $ is not bounded, $X$ has a
complemented subspace isomorphic to 
$\bigl( \oplus _{r=1}^\infty 
(\oplus _{k=1}^r \ell _p^r)_{\ell _q}\bigr) _{\ell _q} .$
Hence we have for each $r$ that $P_r = U_r V_r,$ where
$U_r$ is in ${\Cal F}((\oplus _{k=1}^\infty\ell _p )_{\ell
_q}, X)$ and $V_r$ in ${\Cal F}(X, 
(\oplus _{k=1}^\infty\ell _p )_{\ell _q}),$ $||U_r || = 1 =
||V_r||$ and $U_rV_s =0$ if $r\ne s.$ It follows that $T$
factors as required. 

Now ${\Cal F}( (\oplus _{k=1}^\infty\ell _p)_{\ell _q})$
is amenable, see Corollary 4.8. The above remarks about
 factoring approximable operators therefore allow
us to apply Theorem 6.2 to prove

\proclaim{Theorem 6.5 } Let $\{ n_k\} _{k=1}^\infty $ be 
a sequence of positive integers. Then the algebra ${\Cal F}(
(\oplus _{k=1}^\infty \ell _p^{n_k} )_{\ell _q})$ is
amenable.
\endproclaim

 It was remarked above that 
${\Cal F}(\ell_p\oplus \ell_2)$ is amenable. This suggests
that ${\Cal F}(\ell_p\oplus \ell_q)$ may be amenable for
all $p$ and $q.$ That this is not so will follow
from a further general result about amenable Banach
algebras.

\proclaim{ Definition 6.6} A Banach algebra ${\Cal B}$
has trivial virtual centre if, for each $b^{**}$ in 
${\Cal B}^{**}$ with $bb^{**}=b^{**}b$ for all $b$ in
${\Cal B},$ there is $\lambda $ in ${\Bbb C}$ with
$bb^{**}=\lambda b =b^{**}b$  for all $b$ in ${\Cal B}.$
\endproclaim

The algebras in which we are interested have this property.

\proclaim{Proposition 6.7} Let $X$ be a Banach space.
Then ${\Cal F}(X)$ has trivial virtual centre.
\endproclaim

\demo{ Proof } Let $P$ be a rank one projection on $X.$
Then $P{\Cal F}(X)P={\Bbb C}P.$ Suppose that $B^{**}_0$ in
${\Cal F}(X)^{**}$ satisfies $BB^{**}_0=B^{**}_0B$ for all
$B$ in ${\Cal F}(X).$ Then
$PB^{**}_0=P^2B^{**}_0=PB^{**}_0P.$ Since the map
$B^{**}\mapsto PB^{**}P$ is the second adjoint of the map
$B\mapsto PBP,$ it follows that there is $\lambda_0 $ in
${\Bbb C}$ such that $PB^{**}_0=\lambda_0 P.$ Consequently $\{ T\in
{\Cal F}(X)\mid TB^{**}_0 = \lambda_0 T\} $ is a non-zero 
closed two-sided ideal in the simple Banach algebra ${\Cal
F}(X).$ Therefore $TB^{**}_0 = \lambda_0 T$ for all $T$ in 
${\Cal F}(X).$

\enddemo

For the next theorem let ${\Cal A}$ and ${\Cal A}_{ii}$
be as above.

\proclaim{Theorem 6.8}  Suppose that ${\Cal A}$ is
amenable, that ${\Cal A}_{11}$ and ${\Cal A}_{22}$ have
trivial virtual centre and that ${\Cal A}_{21}$ and 
${\Cal A}_{12}$ are not both zero. Then ${\Cal A}_{jj} =
{\Cal A}_{jj}^\circ $ for at least one value of $j.$ 
\endproclaim

\demo{ Proof } Denote ${\Cal A}^\circ = \{ a\in {\Cal
A}\mid  P_iaP_i \in {\Cal A}_{ii}^\circ ,\ i=1,2 \} .$  On
${\Cal A}^\circ $ define the norm $\| a\| ^\circ = \hbox{
max}\{ \| P_1aP_1\| ^\circ ,\| P_1aP_2 \| ,\| P_2aP_1\|  ,
\| P_2aP_2 \| ^\circ \} .$ Then, for $a\in {\Cal A},\
a^\circ \in {\Cal A}^\circ $ we have $\| aa^\circ \| ^\circ \leq 2\|
a\| \| a^\circ \| ^\circ $ and $\| a^\circ a\| ^\circ \leq 2\| a\| \|
a^\circ \| ^\circ .$ Hence $({\Cal A}^\circ ,\| .\| ^\circ )$ is a
Banach ${\Cal A}$-bimodule.

\par
The map $a\mapsto P_1aP_2 - P_2aP_1 = P_1a - aP_1$ is a derivation
from ${\Cal A}$ into ${\Cal A}^\circ $ and so there is $C$ in $({\Cal
A}^\circ)^{**}$ such that $P_1aP_2 - P_2aP_1 = aC-Ca$ for all $a$ in
${\Cal A}.$ Since ${\Cal A}^\circ = \oplus _{i,j=1,2}P_i {\Cal
A}^\circ P_j ,$ we have $C=\sum _{i,j=1,2} C_{ij},$ where $C_{ij} $
belongs to $(P_i{\Cal A}^\circ P_j)^{**}.$
\par
If $a_{ii}$ belongs to ${\Cal A}_{ii},$ then
$a_{ii}C-Ca_{ii}=0.$ In particular, 
$a_{ii}C_{ii}-C_{ii}a_{ii}=0$ for each $a_{ii}$ in ${\Cal
A}_{ii},$ where $C_{ii}$ belongs to $({\Cal A}_{ii}^\circ
)^{**}.$ Now the second adjoint of the inclusion map ${\Cal
A}_{ii}^\circ \to {\Cal A}_{ii}$ embeds $({\Cal
A}_{ii}^\circ )^{**}$ in $({\Cal A}_{ii})^{**}$ and so,
since ${\Cal A}_{ii}$ has trivial virtual centre for each
$i,$ there are $\lambda _1$ and $\lambda _2$ in ${\Bbb C}$
such that $a_{ii}C_{ii}=\lambda _ia_{ii} =C_{ii}a_{ii}$ for
$a_{ii}$ in ${\Cal A}_{ii},$ $i=1,2.$ 
\par
Suppose, without loss of generality, that ${\Cal A}_{12}$
is not zero and choose $a_{12}\ne 0$ in ${\Cal A}_{12}.$ 
Since ${\Cal A}$ is amenable, it has a bounded two-sided
approximate identity and so, by Proposition 6.1, ${\Cal
A}_{12}$ is an essential left ${\Cal A}_{11}$- and
essential right ${\Cal A}_{22}$-module. Hence there are
$a_{11}$ in ${\Cal A}_{11},$ $a_{22}$ in ${\Cal A}_{22}$
and $\tilde a_{12}$ in ${\Cal A}_{12}$ with
$a_{12}=a_{11}\tilde a_{12}a_{22}.$ We have
$a_{12}=a_{12}C-Ca_{12}=a_{12}C_{22}-C_{11}a_{12}.$
Substituting for $a_{12}$ we get $a_{12}=
a_{11}\tilde a_{12}a_{22}C_{22}-C_{11}
a_{11}\tilde a_{12}a_{22} = \lambda _2
a_{11}\tilde a_{12}a_{22}-\lambda _1
a_{11}\tilde a_{12}a_{22} =(\lambda _2-\lambda _1)a_{12}.$
Therefore $\lambda _2-\lambda _1 =1$ and so at least one
of $\lambda _1$ and $\lambda _2$ is not zero.
\par

Suppose that $\lambda _1\ne 0$ and put $b^{**}=\lambda
_1^{-1}C_{11}.$ Then $b^{**}$ belongs to $({\Cal
A}_{11}^\circ )^{**}$ and $a_{11}b^{**}=a_{11}$ for every
$a_{11}$ in ${\Cal A}_{11}.$ Hence, if $\{ b^\alpha \}
_{\alpha \in A}$ is a bounded net in ${\Cal A}_{11}^\circ
$ which converges to $b^{**}$ in the $\text{weak}^*$- topology, then
$\{ a_{11}b^\alpha \} _{\alpha \in A}$ converges weakly to
$a_{11}.$ It follows, as in \cite{B\&D,Proposition 11.4},
that there is a net $\{ e^\beta \} _{\beta \in B},$ each
$e^\beta $ being a convex combination of $b^\alpha $'s,
which is a right approximate identity for ${\Cal
A}_{11}.$ The approximate identity $\{ e^\beta \} _{\beta
\in B}$ is bounded, by $\| b^{**}\| ,$ in $({\Cal
A}_{11}^\circ , \| .\| ^\circ )$ and ${\Cal
A}_{11}^\circ $ is an ideal in ${\Cal A}_{11}.$ Hence for
each $a_{11}$ in ${\Cal A}_{11}$ and each $\epsilon >0$
there is $c=a_{11}e^\beta $ with $\| a_{11}-c\| <\epsilon
$ and $\| c\| ^\circ <2\| b^{**}\| \| a_{11}\| .$ Consequently,
for each $a_{11}$ in ${\Cal A}_{11},$ there is a series
$\sum_i c_i $ in ${\Cal A}_{11}^\circ $ with $\sum _i\|
c_i\| ^\circ <\infty $ and $\sum_i c_i = a_{11}.$
Therefore ${\Cal A}_{11}^\circ ={\Cal A}_{11}.$

\enddemo

These last two results may be reformulated to say that
the spaces $X$ with ${\Cal F}(X)$ amenable have a
property which is a little like being primary. Recall
that a Banach space, $X,$ is {\it primary } if, for every
bounded projection $Q$ on $X,$ either $QX$ or $(I-Q)X$ is
isomorphic to $X,$ see \cite{L\& T, Definition 3.b.7}.
Let us say that $X$ is {\it approximately primary } if,
for every bounded projection $Q$ on $X,$ at least one of
the product maps $\pi :  {\Cal F}(X, QX)\pten {\Cal
F}(QX, X) \to {\Cal F}(X) $ or $\pi : {\Cal F}(X,(I-
Q)X)\pten {\Cal F}((I-Q)X, X) \to {\Cal F}(X) $ is
surjective. Then every primary space is approximately
primary as is every space with a subsymmetric basis, see
\cite{L\&T, Proposition 3.b.8}. 

 Now put ${\Cal A}={\Cal F}(X)$ and suppose
that ${\Cal A}$ is amenable. Let $P_1$ be the idempotent
in $M({\Cal A})$ determined by a bounded projection
$Q$ on $X.$ Then ${\Cal A}_{11}$ is
isomorphic to ${\Cal F}(QX)$ and  ${\Cal A}_{22}$ to ${\Cal
F}((I-Q)X).$ Hence, by Proposition 6.7,  ${\Cal A}_{ii}$
has trivial virtual centre for $i=1,2.$ Clearly, ${\Cal
A}_{12}$ is not zero and so, by Theorem 6.8, ${\Cal
A}_{ii}^\circ ={\Cal A}_{ii}$ for at least one value of
$i.$  It follows that, if ${\Cal F}(X)$ is amenable, then
$X$ is approximately primary.

\proclaim{Theorem 6.9} If $1<p,q<\infty ,$ $p\ne q$ and neither
$p$ nor $q$ is equal to 2, then ${\Cal F}(\ell_p\oplus \ell_q)$
is not amenable.
\endproclaim

\demo{ Proof } In view of 6.7, 6.8 and the remarks
following it suffices to show that, if $1<p,q<\infty ,$ $p
\ne q$ and neither is equal to 2, then the product map
$\pi : {\Cal F}(\ell_p, \ell _q)\pten {\Cal F}(\ell
_q, \ell _p) \to {\Cal F}(\ell _p) $ is not surjective. 

Suppose that $\pi $ is surjective. Then, by the open
mapping theorem, there is a $K>0$ such that for each
$T$ in ${\Cal F}(\ell _p)$ we have $T= \pi (\sum
_{n=1}^\infty U_n \otimes V_n ),$ where
$\sum _{n=1}^\infty ||U_n|||| V_n || < K||T||.$ It
follows, since $\ell _q$ is isomorphic to $(\oplus
_{n=1}^\infty \ell _q)_{\ell _q},$ that $T=UV,$ where $U$
is in ${\Cal F}(\ell_p, \ell _q),$ $V$ in 
${\Cal F}(\ell _q, \ell _p)$ and $||U||||V|| < K||T||.$ 

Let $P_j$ be the projection onto the span of the first
$j$ vectors of the standard basis for $\ell _p.$ Then,
since we are supposing that $\pi $ is surjective, $P_j =
U_j V_j $ where $||U_j||||V_j|| < K.$ Put $Q_j =
V_jU_j.$ Then $Q_j$ is a projection on $\ell _q$ and
$||Q_j|| < K.$ Defining $U_j'=P_jU_jQ_j$ and
$V_j'=Q_jV_jP_j,$ we have that $U_j'$ is an
isomorphism from the range of $Q_j$ to the range of
$P_j,$ $V_j'$ is the inverse of $U_j'$ and
$||U_j'||||V_j'|| < K^3.$ Hence, if $\pi $ is surjective,
then $\ell _p$ is finitely representable in $\ell _q,$ see
\cite{Wo, Definition II.E.15}. It is known that this is not so
if $p\ne q$ and neither is equal to 2. There are several
cases. 

First, suppose that $p<2<q.$ If $\ell _p$ were finitely 
representable in $\ell _q,$ then, since $\ell _q$ is of
type 2, $\ell _p$ would be of type 2. (See \cite{Wo, 
Definition III.A.17 and Theorem III.A.23} ) That is not
so. Therefore $\ell _p$ is not finitely representable in
$\ell _q.$ The case $q<2<p$ is dual to this case. 

Next, suppose that $2<q<p.$ If $\ell _p$ were finitely 
representable in $\ell _q,$ then, since $\ell _q$ is of
cotype q, $\ell _p$ would be of cotype q. Since $\ell
_p$ is not of cotype q, $\ell _p$ is not finitely
representable in $\ell _q.$ The case $p<q<2$ is dual to
this case.

Finally, suppose that $2<p<q.$ If $\ell _p$ were finitely
representable in $\ell _q,$ then $\ell _p$ would be isomorphic to a
subspace of $L_q(\mu )$ for some measure $\mu ,$ see Proposition 7.1
in [L\&P]. It would then follow, by Corollary 2 in [K\&P], that $\ell
_p$ had a complemented subspace isomorphic to $\ell _2$ or $\ell _q.$
However, that is not possible because, by Proposition 2.c.3 in [L\&T],
every operator from $\ell _p$ to $\ell _2$ and every operator from
$\ell _q$ to $\ell _p$ is compact. The case $q<p<2$ is dual to this
case. This argument is also sketched on [Wo] pages 104 and 107.

\enddemo

The above proof also shows that ${\Cal F}(c_0\oplus \ell
_p)$ is not amenable when $p<2$ but does not treat the
case $p>2.$  Similarly, ${\Cal F}(\ell
_1\oplus \ell _p)$ is not amenable when $p>2.$ 

We conclude this section with a result which shows that
amenability of ${\Cal F}(X)$ is partially preserved on
complemented subspaces of $X.$ 

\proclaim{ 6.10 Theorem } Let $X$ and $Y$ be Banach spaces
and suppose that ${\Cal F}(X\oplus Y)$ is amenable. Then
at least one of ${\Cal F}(X)$ and ${\Cal F}(Y)$ is
amenable.
\endproclaim

\demo{ Proof } Put ${\Cal A}={\Cal F}(X\oplus Y)$ and let
$P_1$ be the idempotent in $M({\Cal A})$ determined by
the projection onto $X$ with kernel $Y.$ Then, by 6.7 and
6.8, ${\Cal A}_{jj} = {\Cal A}_{jj}^\circ $ for at least
one value of $j.$ By Theorem 6.2, it follows that at
least one of ${\Cal A}_{11}$ and ${\Cal A}_{22}$ is
amenable. Since ${\Cal A}_{11}$ is isomorphic to ${\Cal
F}(X)$ and ${\Cal A}_{22}$ is isomorphic to ${\Cal
F}(Y),$ the result follows.

\enddemo

The conclusion of this last theorem is the best possible,
that is, there are spaces $X$ and $Y$ such that ${\Cal
F}(X\oplus Y)$ is amenable but ${\Cal F}(X)$ is not. For
example, let $X=c_0\oplus \ell_1$ and $Y=\ell_1(c_0)= 
\bigl( \oplus _{n=1}^\infty c_0\bigr) _{\ell_1}.$ Then
$X\oplus Y$ is isomorphic to $Y.$ Hence ${\Cal F}(X\oplus Y)$ and
${\Cal F}(Y)$ are amenable by Corollary 4.8. On the other hand, ${\Cal
F}(X)$ is not amenable by Theorem 6.9.

\heading
{7. Open questions and conclusion }
\endheading

The name `amenable' is used for a Banach algebra $\Cal A$ satisfying
the cohomological condition $H^1(\Cal A,X^*)=0$ for all Banach- $\Cal
A$-modules $X$, see [J], because of the theorem that a group algebra
$L^1(G)$ satisfies this condition if and only if the locally compact
group $G$ is amenable, [J]. Amenability is an important property of
groups which has many characterizations. As well as the cohomological
characterization of the group algebra, it may be described in terms of
group representations, fixed points of group actions, translation
invariant functionals and in other ways.  The F{\o}lner conditions on
compact subsets of the group characterize amenability in terms of
properties intrinsic to the group. Alternative characterizations of
the amenability of ${\Cal K}(X)$ and ${\Cal F}(X)$ would help us to
have a better understanding of its significance. We are thus led to
ask

\proclaim {Question 7.1} What are the intrinsic properties
of the Banach space $X$ which are equivalent to amenability of ${\Cal
K}(X)$ and ${\Cal F}(X)?$
\endproclaim

The results we have obtained so far suggest that
amenability of  ${\Cal K}(X)$ and ${\Cal F}(X)$ may be
equivalent to some sort of approximation property for
$X.$ Such an approximation property, if it were to exist,
would be the analogue of the F\o lner conditions. 
\par

Approximation properties are certainly necessary. Since an amenable
algebra has a bounded two-sided approximate identity, if ${\Cal K}(X)$
is amenable, then $X^*$ has what is called in [G\&W] the $B -{\Cal
K}(X)^a $- AP, and in [Sa] the $*$ - b.c.a.p., that is, the identity
operator on $X^*$ is approximable in the topology of convergence on
compacta by operators which are adjoints of compact operators on $X$.
It also follows, by \cite{D, Theorem 2.6}, that $X$ has the bounded compact
approximation property.  Similarly, if ${\Cal F}(X)$ is amenable, then
$X$ and $X^*$ have the bounded approximation property. However, the
relationship between amenability of ${\Cal K}(X)$ and the
approximation property is not clear.

\proclaim{Question 7.2} Does amenability of ${\Cal K}(X)$
imply that $X$ has the approximation property?
\endproclaim

\noindent If there should be a Banach space $X$ which does
not have the approximation property but is such that ${\Cal K}(X)$ is
amenable , then ${\Cal K}(X)/ {\Cal F}(X)$ would be a radical,
amenable Banach algebra. At present no example of such a Banach
algebra is known.

\par
Theorem 6.9  shows that for ${\Cal F}(X)$ to be amenable
it does not suffice that $X^*$ have the $B -{\Cal
F}(X)^a $- AP. It seems necessary for there also to be
some sort of symmetrization of the approximation property.
We have seen, in Section 4, a symmetrized approximation
property, property (\A), which forces the amenability of ${\Cal F}(X).$
This property was used to show  that ${\Cal F}(X)$ is
amenable for many of the classical Banach spaces and for
spaces with a shrinking, subsymmetric basis. 

\proclaim{Question 7.3} Is property (\A) or some similar symmetrized
approximation property equivalent to amenability of ${\Cal F}(X)$ or 
${\Cal K}(X)?$
\endproclaim

In order to determine how close this property is to being equivalent
to the amenability of ${\Cal F}(X),$ it would be useful to investigate
whether ${\Cal F}(X)$ is amenable if $X$ is a space which is clearly
unlikely to have this symmetrized approximation property. Examples
that come to mind are the James space, which does not have an
unconditional basis (\cite{L\&T, 1.d.2}), and the Tsirelson space,
which contains no subsymmetric basic sequence (\cite{L\&T, p. 132}).

\proclaim{Question 7.4} Is ${\Cal F}(X)$ amenable if $X$
is the James space or the Tsirelson space?
\endproclaim

We have seen that the class of spaces, $X,$ such that 
${\Cal F}(X)$ is amenable is not closed under direct sums
or under passing to complemented subspaces. However, any
space, $X,$ such that ${\Cal F}(X)$ is amenable has the
property that $X^*$ satisfies the $B -{\Cal F}(X)^a
$- AP and this property is inherited by complemented
subspaces of $X$ and is preserved under direct sums.
Perhaps this is the most that can be said about such
spaces.

\proclaim{Question 7.5} Is the smallest space ideal 
containing all spaces, $X,$ such that ${\Cal F}(X)$ is
amenable equal to the class of all Banach spaces, $X,$
such that $X^*$  has the  $B -{\Cal F}(X)^a $- AP?
\endproclaim

\noindent Recall from \cite{Pie, Definition 2.1.1}, that a
{\it space ideal } is a class of Banach spaces which
contains the finite dimensional spaces and is closed under
direct sums and taking complemented subspaces. It is
clear that the class of spaces such that $X^*$ has the 
$B -{\Cal F}(X)^a $- AP is a space ideal. Should
the answer to 7.2 be `no', an obvious further question
would be whether the class of all Banach spaces whose
duals have the  $B -{\Cal K}(X)^a $- AP is equal to
the smallest space ideal containing all spaces, $X,$ such
that ${\Cal K}(X)$ is amenable.

It was shown by J. Lindenstrauss, see \cite{L\&T, Theorem
3.b.1}, that every Banach space with an unconditional basis
is isomorphic to a complemented subspace of a space with a
symmetric basis. In view of the possible equivalence of
the amenability of ${\Cal F}(X)$ with some symmetric
approximation property, this suggests the following
refinement of 7.5.

\proclaim{ Question 7.6 } Is every Banach space, $X,$
such that $X^*$ has the  $B -{\Cal F}(X)^a $- AP 
isomorphic to a complemented subspace of a space, $Y,$
such that  ${\Cal F}(Y)$ is amenable?
\endproclaim

The spaces $C_p,$ $1\leq p<\infty,$ introduced by W. B.  Johnson
\cite{Jo1} will provide the answer to this last question. The space
$C_p$ is the $\ell ^p$ direct sum of a sequence of finite dimensional
spaces which is dense in the set of all finite dimensional spaces. It
has the property that every approximable operator factors through it
and, for $1<p<\infty,$ ${\Cal F}(C_p)$ has a bounded two-sided
approximate identity. Now let $X$ be any space such that $X^*$ has the
$B -{\Cal F}(X)^a $- AP. Then, since $C_p$ has the above properties,
Theorem 5.2 implies that ${\Cal F}(X\oplus C_p)$ is amenable if and
only if ${\Cal F}(C_p)$ is amenable. Therefore the answer to 7.6 is `yes' if
${\Cal F}(C_p)$ is amenable. On the other hand, if $C_p$ is isomorphic
to a complemented subspace of some space, $Y,$ such that ${\Cal F}(Y)$
is amenable, that is, if the answer to 7.6 is `yes' when $X = C_p,$
then ${\Cal F}(C_p)$ is amenable.

\proclaim{ Question 7.7 } Is ${\Cal F}(C_p)$ amenable for
any, and hence all, $1<p<\infty ?$ 
\endproclaim

\noindent Note that $C_1^*$ does not have the approximation property,
see \cite{Jo2, Theorem 3}, and so ${\Cal F}(C_1)$ is not amenable.

Another theorem, similar to that of Lindenstrauss, is
proved in [J,R\&Z] and [P], see \cite{L\&T, Theorem 1.e.13}. 
It says that any separable Banach space with the
B.A.P. is isomorphic to a complemented subspace of a
Banach space with a basis. There is an even stronger
theorem, see \cite{L\&T, Theorems 2.d.8 and 2.d.10}, that there
is a Banach space, $U,$ with basis such that any separable
Banach space with the B.A.P. is isomorphic to a
complemented subspace of $U$ and that $U$ is determined
uniquely up to isomorphism by this property. The space
$U$ is said to be complementably universal for the spaces
with the B.A.P. Now if $X$ has a shrinking basis, then
$X^*$ has the $B -{\Cal F}(X)^a $- AP. This suggests 

\proclaim{ Question 7.8 } (a) Is there a Banach space, $V,$
with a shrinking basis which is complementably universal 
for the spaces, $X,$ such that $X^*$ has
the $B -{\Cal F}(X)^a $- AP? \hfill \break
(b) If so, is ${\Cal F}(V)$ amenable?
\endproclaim
\medskip
\subheading{Acknowledgements} We are indebted to D\. J\. H\. Garling for
his help with a proof of Theorem 6.9 and to P\. G\. Casazza, J\.
Partington and A\. Pe\l czy\'nski for valuable discussions. We would
also like to thank  L\. Tzafriri for having brought the reference
\cite{A\&F} to our attention.  

This work was initiated while the first and the third author were
visiting University of Leeds. We are both thankful for the hospitality
extended.

\enddocument